\documentclass[reqno,11pt]{amsart}
\usepackage{amsmath,amsthm,amscd,amsfonts,amssymb,hyperref,floatrow}
\usepackage{graphicx, color,dsfont,xcolor}

\numberwithin{equation}{section}

\newtheorem{theorem}{Theorem}[section]
\newtheorem{lemma}[theorem]{Lemma}
\newtheorem{proposition}[theorem]{Proposition}
\newtheorem{corollary}[theorem]{Corollary}
\newtheorem{rem}[theorem]{Remark}

\renewcommand{\tilde}{\widetilde}          
\DeclareMathSymbol{\leqslant}{\mathalpha}{AMSa}{"36} 
\DeclareMathSymbol{\geqslant}{\mathalpha}{AMSa}{"3E} 
\DeclareMathSymbol{\eset}{\mathalpha}{AMSb}{"3F}     
\renewcommand{\leq}{\;\leqslant\;}                   
\renewcommand{\geq}{\;\geqslant\;}                   

\newcommand{\R}{\mathbb{R}}

\newcommand{\N}{\mathbb{N}}

\def \P{ \mathbb P  }
 \newcommand{\eps}{\varepsilon}
\def \E{ \mathbb E  }

\begin{document}
\title{From Sine kernel to Poisson statistics}

\author{Romain Allez \and Laure Dumaz}
\address{Weierstrass Institute, Mohrenstr. 39, 10117 Berlin, Germany.}
\address{Statistical Laboratory, Centre for Mathematical Sciences, Wilberforce Road, Cambridge, CB3 0WB, United Kingdom. }

\email{romain.allez@gmail.com,L.dumaz@statslab.cam.ac.uk} 
\date{\today}

\maketitle
\begin{abstract}
We study the Sine$_\beta$ process 
 introduced in [B. Valk\'o and B. Vir\'ag.  
Invent. math. {\bf 177} 463-508 (2009)] when the inverse temperature $\beta$ tends to $0$. 
This point process has been shown to be the scaling limit of the eigenvalues point process 
in the bulk of $\beta$-ensembles and its law is characterized in terms of 
the winding numbers of the Brownian carrousel at different angular speeds.   
After a careful analysis of this family of coupled diffusion processes, we prove that the Sine$_\beta$ point process 
converges weakly to a Poisson point process on $\R$.
Thus, the Sine$_\beta$ 
point processes establish a smooth crossover between the rigid clock (or picket fence) process (corresponding to $\beta=\infty$) and the Poisson process. 
\end{abstract}
\vspace{1cm}

\section{Introduction and main result} 
Although random matrices were originally introduced by John Wishart \cite{wishart} in 1928 as a tool to study 
population dynamics in biology through principal component analysis, they became very popular 
much later in 1951 when Wigner \cite{wigner} postulated that the fluctuations in positions of the energy levels of heavy nuclei are well described
(in terms of statistical properties) by the eigenvalues of a very large Hermitian random matrix.  
Random matrix theory (RMT) is now an active research area in mathematics and theoretical physics with applications in 
statistics, biology, financial mathematics, engineering and telecommunications, number theory etc. (see \cite{agz,silverstein,mehta,forrester,handbook} for a state of the art).

The classical models of Hermitian random matrices are the Gaussian orthogonal, unitary and symplectic ensembles. 
It is well known that the joint law of the eigenvalues of the matrices in those Gaussian ensembles is the 
Boltzmann-Gibbs equilibrium measure of a one-dimensional repulsive Coulomb gas confined in a harmonic well. 
More precisely, this joint law has a probability density $\mathcal{P}_\beta$ on $\R^N$ ($N$ is the dimension of the square matrices) 
given by 
\begin{align}\label{def-Pbeta}
\mathcal{P}_\beta(\lambda_1,\cdots,\lambda_N) = \frac{1}{Z_N} \prod_{i <j} |\lambda_i-\lambda_j| ^{\beta} \exp(-\frac{N\beta}{4}\sum_{i=1}^N \lambda_i^2)
\end{align}
where the inverse temperature $\beta=1$ for the Gaussian orthogonal ensemble, respectively $\beta=2,4$ for
the unitary and symplectic ensembles. 
The linear statistics of the point processes with joint probability density functions (jpdf) $\mathcal{P}_\beta,\beta=1,2,4$ have been 
extensively studied in the literature with different methods \cite{agz,mehta}. 

In 2002, Dumitriu and Edelman \cite{dumitriu} came up with a new explicit ensemble of random tri-diagonal matrices whose eigenvalues
are distributed according to the jpdf $\mathcal{P}_\beta$ for any $\beta>0$ (see also \cite{alice,jp-alice} where invariant
ensembles associated to general $\beta$ were constructed). 

Those tri-diagonal matrices have been very useful in the last decade, leading to important progress 
on the understanding of the local eigenvalues statistics in the limit of large dimension $N$ for general $\beta>0$. 
At the edge of the spectrum, 
it was first proved \cite{virag-1} that the largest eigenvalues converge jointly 
(when zooming in the edge-scaling region of width $N^{-2/3}$ around $2$)
to the low lying eigenvalues of a random Schr\"odinger operator called the stochastic Airy operator (see also \cite{sutton}).  
Similar results were proved for the bulk in \cite{virag} by Valk\'o and Vir\'ag. 
For $\lambda$ belonging to the Wigner sea $(-2,2)$, the authors of \cite{virag} consider the point process 
\begin{align}\label{eigenvalues-rescaled-pp}
\Lambda_N:= \left( 2\pi N \rho(\lambda) \left(\lambda_i - \lambda \right) \right)_{i=1,\cdots,N}
\end{align}
where $ (\lambda_1,\cdots, \lambda_N)$ is distributed according to $\mathcal{P}_\beta$ and 
$\rho(\lambda)= \frac{1}{2\pi} \sqrt{4-\lambda^2}$ is the Wigner semi-circle density.   
Indeed, the mean level spacing around level $\lambda$ for the points with law $\mathcal{P}_\beta$ is approximately 
$1/(N\rho(\lambda))$ when $N \gg 1$. The mean point spacing 
of $\Lambda_N$ defined in \eqref{eigenvalues-rescaled-pp} is therefore of order $2\pi$ and in this scaling, one can now investigate 
the limiting statistics of this point process when $N\to\infty$. The authors of \cite{virag} precisely answer this question 
proving that the point process $\Lambda_N$ converges in law 
\footnote{The convergence is with respect to vague topology for the counting measure of the point process.}
to a point process Sine$_\beta$ on $\R$ first introduced in \cite{virag} and 
characterized in terms of a family $(\alpha_\lambda)_{\lambda\in \R}$ of coupled 
one-dimensional diffusion processes, the \emph{stochastic sine equations}.  As expected for the eigenvalues statistics in the bulk, 
the point process Sine$_\beta$ is translation-invariant in law. 
The family of diffusions $(\alpha_\lambda)_{\lambda\in \R}$ can be interpreted as the hyperbolic angle of the Brownian carousel 
with parameter $\lambda$ and its law is characterized as follows: 
Given a (driving) complex Brownian motion $(Z_t)_{t\ge 0}$, the diffusions $\alpha_\lambda, \lambda\in \R$ 
satisfy
\begin{align}\label{eq-alpha}
d\alpha_\lambda= \lambda\frac{\beta}{4} e^{-\frac{\beta}{4}t} \, dt + { \rm Re}((e^{-i \alpha_\lambda}-1) dZ_t)\,,  \quad
\alpha_\lambda(0)=0\,.  
\end{align}  
Note that all the diffusions  $\alpha_\lambda, \lambda\in \R$ are adapted to the filtration 
of the Brownian motion $(Z_t)$. This coupling induces a strong interaction between the diffusions which makes 
the joint law difficult to analyse, as we shall see. A key feature shared by the processes 
 $\alpha_\lambda,\lambda\in \R$ is that they all converge almost surely as $t\to\infty$ to a limit $\alpha_\lambda(\infty)$ which is 
an integer multiple of $2\pi$.  
The characterization of the law of the Sine$_\beta$ point process \footnote{If $A$ is a Borel set of $\R$, Sine$_\beta(A)$ inside $A$ denotes the number of points inside $A$. In other words, Sine$_\beta$ is the counting measure of the point process.} can now be enunciated as follows
\begin{align}\label{def-Sine-beta}
\left({\rm Sine}_\beta([\lambda,\lambda'])\right)_{\lambda< \lambda'} \stackrel{(d)}{=} \left( \frac{\alpha_{\lambda'}(\infty)-\alpha_\lambda(\infty)}{2\pi} \right)_{\lambda<\lambda'} . 
\end{align}

\begin{figure}[h!btp] 
     \center
     \includegraphics[width=12cm]{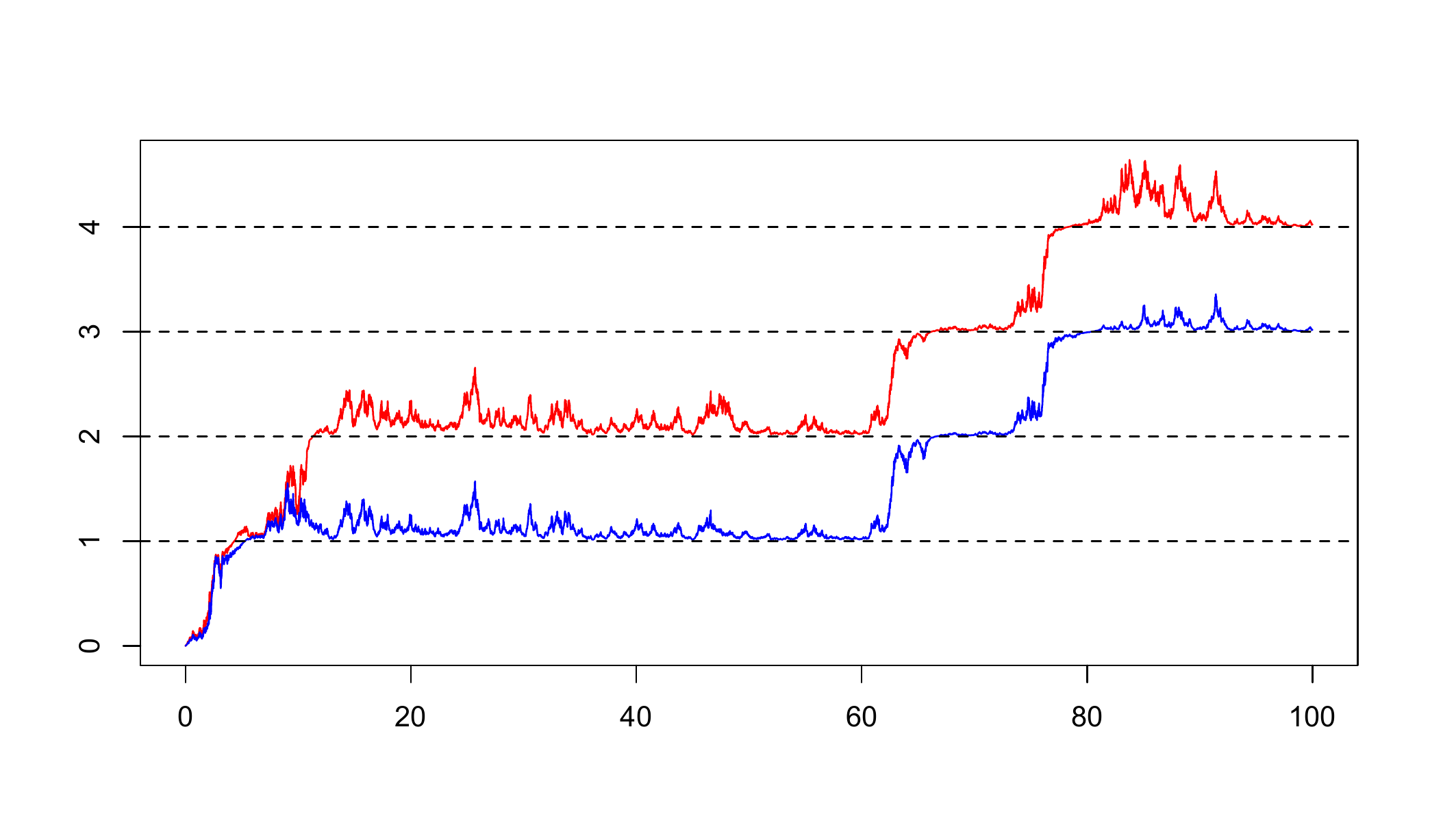}
     \caption{(Color online). Sample paths of two diffusions $\alpha_\lambda/(2\pi)$ (blue curve) and $\alpha_{\lambda'}/(2\pi)$ (red curve) for $\lambda <\lambda'$ and for a small value of $\beta$. 
     We see that whenever $\lfloor \alpha_\lambda/(2\pi) \rfloor$ jumps,  $\lfloor \alpha_{\lambda'}/(2\pi) \rfloor$ jumps at the same time in agreement with Lemma \ref{inclusion-croissante}. }
     \label{without-diff}
\end{figure}

In this paper, we are interested in the limiting law of the Sine$_\beta$ process when the inverse temperature $\beta$ 
goes to $0$. Theorem \ref{main} is the main result of this paper and gives the convergence as $\beta \to 0$
of the Sine$_\beta$ process towards a Poisson point process on $\R$. 
This convergence at the continuous ($N=\infty$) level  seems rather natural since taking $\beta \to 0$ amounts to decreasing 
the electrostatic repulsion (and hence the correlation) at the discrete level, i.e. between the $N$ points distributed according to $\mathcal{P}_\beta$. 
If one takes $\beta= 0$ abruptly for fixed $N$, 
the probability density $\mathcal{P}_0$ corresponds (up to a rescaling depending on $\beta$) to the joint law of independent 
Gaussian variables and it is then straightforward to check the convergence (as $N\to \infty$) 
of the point process with law $\mathcal{P}_0$ towards a Poisson process. 
Theorem \ref{main} exchanges the order of the limits $\beta\to 0$ and $N\to \infty$, describing 
the statistics 
when $N\to \infty$ {\it first} and then $\beta \to 0$. It also gives the precise rate of the convergence. 

\begin{theorem}\label{main}
As $\beta\to 0$, the Sine$_\beta$ point process converges weakly in the space of Radon measure 
(equipped with the topology of vague convergence  \cite{kallenberg}) to a Poisson point process on $\R$ 
with intensity $\frac{d\lambda}{2\pi}  $. 
In particular, we have, for any $k\in \N$ and $\lambda < \lambda'$, 
\begin{align*}
\P[{\rm Sine}_\beta [\lambda,\lambda']=k] \to_{\beta\to 0} \exp\left(\frac{\lambda-\lambda'}{2\pi} \right) \frac{ (\lambda'-\lambda)^k  }{(2\pi )^k k!}\,, 
\end{align*}
and the numbers of points of Sine$_\beta$  inside two disjoint intervals are asymptotically independent. 
\end{theorem}

let us briefly discuss some implications of Theorem \ref{main} and 
mention a few related questions on the spectral 
statistics of random matrices and random Schr\"odinger operators. 

In \cite{killip}, the authors have shown that the circular $\beta$-ensemble, which was later shown 
to be Sine$_\beta$ in \cite{nakano}, 
interpolates between Poisson and clock distributions on the circle 
(point process with rigid spacings like the numerals on a clock) 
by considering random CMV matrices. Theorem \ref{main} provides a more precise description of  
this interpolating process on the Poisson process side. 

In our study, we are led to examine a time homogeneous 
family of diffusions $(\theta_\lambda)_{\lambda\in \R}$  defined as
\begin{align}\label{eq-alpha}
d\theta_\lambda= \lambda\frac{\beta}{4}  \, dt + { \rm Re}((e^{-i \theta_\lambda}-1) dZ_t)\,,  \quad
\theta_\lambda(0)=0\,.  
\end{align}  
This family of coupled diffusions also appears in \cite{kritchevski} 
to describe the law of the limiting point process of a certain critical 
random discrete Schr\"odinger operator. Our result can be extended in this context to prove that 
this critical Schr\"odinger operator continuously interpolates between the extended 
(clock/picket fence) and localized (Poisson) regimes. More precisely, one could prove 
using our ideas that the random spectrum has Poissonian statistics in the limit of large temperature.

Let us also compare the results stated in Theorem \ref{main} with those of a previous work \cite{laure-romain} where we consider 
the stochastic Airy ensemble, Airy$_\beta$, obtained in the scaling limit of $\beta$-ensembles at the edge of the spectrum. 
In this context, we proved that the number of points 
Airy$_\beta(]-\infty,\lambda] )$ inside the interval $]-\infty,\lambda]$ displays Poisson statistics in the small $\beta$ limit. 
This permitted us to obtain the limiting distributions  as $\beta\to 0$ 
of each of the lowest eigenvalues (individually) of the Airy$_\beta$ ensemble.  
In particular, we obtained the weak convergence of the $TW(\beta)$ distribution towards the Gumbel distribution.
 Although the Sine$_\beta$ and Airy$_\beta$ characterizations in law (in terms of a family of coupled diffusions) 
look very similar, the analysis of the limiting marginal statistics of the number of points inside a finite closed interval
Airy$_\beta[\lambda,\lambda'] $ for $\lambda < \lambda'$ and the asymptotic independence  when $\beta\to 0$
of the respective numbers of points of the Airy$_\beta$ point process into two disjoint intervals remain open even after our study 
\cite{laure-romain}. 
In this aspect, Theorem \ref{main} gives a much more powerful and complete 
description of the Sine$_\beta$ process in the small $\beta$ limit. 
In this case, we are able to prove the asymptotic independence between the number of points of the Sine$_\beta$ process
in two disjoint intervals. This part of the proof requires new ideas in order to obtain estimates on the relative positions between two coupled
diffusions $\alpha_\lambda$ and $\alpha_{\lambda'}$.  
The nice feature of the Sine$_\beta$ process is its translation invariance in law. This property 
makes the analysis of Sine$_\beta$ easier than the one of the Airy$_\beta$ process. 
The non-homogeneous intensity of the Airy$_\beta$ process is governed by the edge-scaling crossover spectral density 
 of $\beta$-ensembles computed explicitly in \cite{bb,forrester.tw} for $\beta=2$ (see also \cite{laure-romain}).

 
Other spectral statistics of random matrices 
at high temperature, i.e. when $\beta \to 0$, have been investigated in \cite{jp-alice,satya}. 
In \cite{jp-alice}, the authors study the empirical eigenvalue density in the limit of large dimension $N$ for $\beta$-ensembles when 
$\beta$ tends to $0$ with $N$ as $\beta= 2c/N$ where $c>0$ is a constant. 
The authors compute the limiting spectral density $\rho_c(\lambda)$ explicitly in terms of parabolic cylinder function 
and establish a Gauss-Wigner crossover, in the sense that the family $\rho_c$ interpolate between the Gaussian probability distribution ($c=0$)
and the Wigner semi-circle ($c\to +\infty$). The case of Gaussian Wishart matrices has also been studied in \cite{satya}. 

It would be interesting to have a description of the crossover statistics of the $\beta$-ensembles 
obtained in the double scaling limits
when $\beta$ tends to $0$ with $N\to \infty$ (this question is briefly discussed in \cite{laure-romain} for the statistics 
at the edge of the spectrum).

\begin{figure}[h!btp] 
     \center
     \includegraphics[width=12cm]{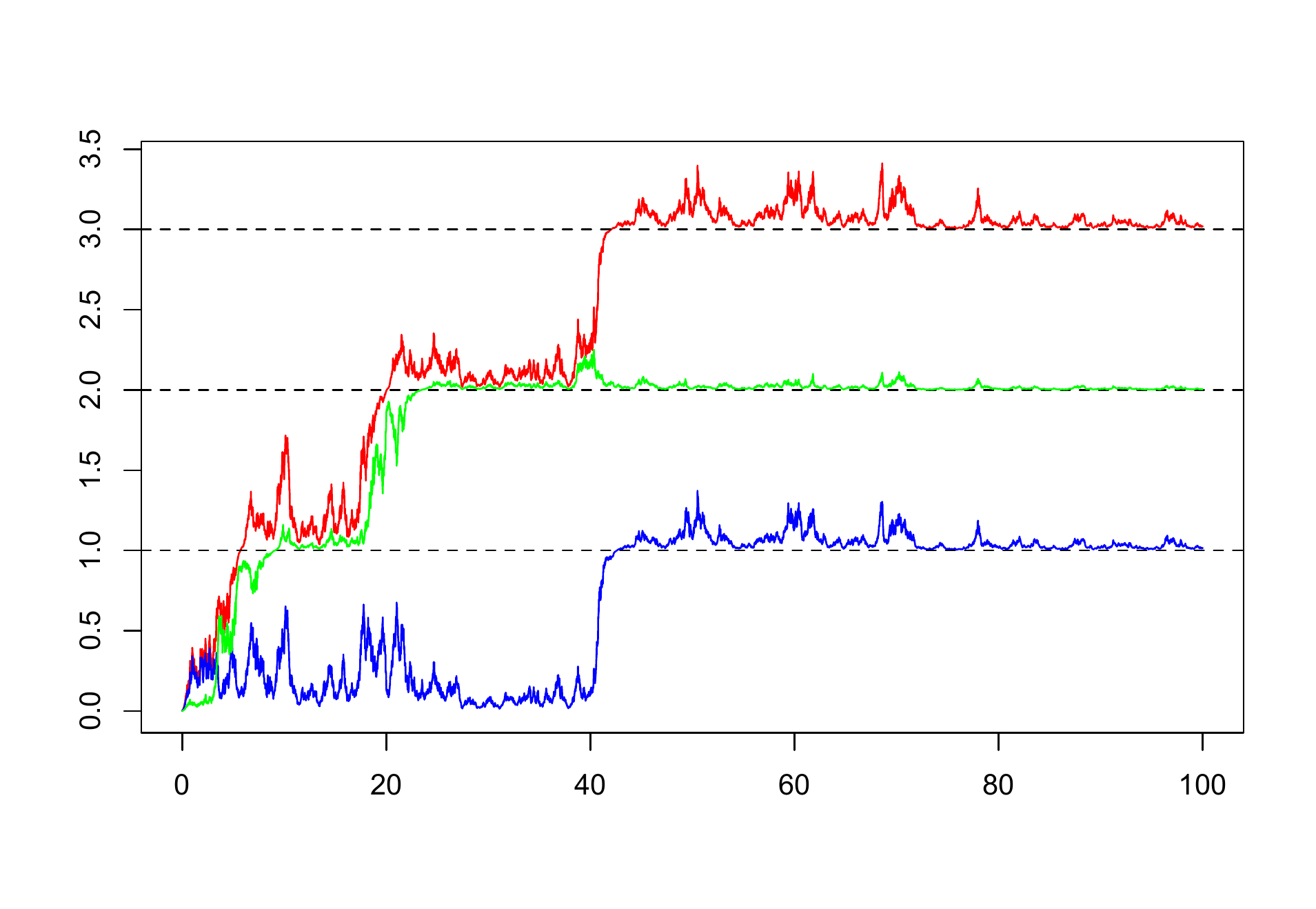}
     \caption{(Color online). Sample paths of the diffusions $\alpha_\lambda/(2\pi)$ (blue curve) and $\alpha_{\lambda'}/(2\pi)$ (red curve) together with $(\alpha_{\lambda'}-\alpha_\lambda)/(2\pi)$
     (green curve)  
     for $\lambda <\lambda'$ and for a small value of $\beta$. 
     Again we see that $\lfloor \alpha_\lambda/(2\pi) \rfloor$ and $\lfloor (\alpha_{\lambda'} -\alpha_\lambda) /(2\pi) \rfloor$ never jump simultaneously while  $\lfloor \alpha_{\lambda'}/(2\pi) \rfloor$ and $\lfloor (\alpha_{\lambda'} -\alpha_\lambda) /(2\pi) \rfloor$ always jump  at the same times in agreement with Lemma 
     \ref{poissons-indep}, Lemma \ref{inclusion-croissante} and  Remark \ref{remark-inclusion-facult}.  } \label{with-diff}
\end{figure}

{\it Organization of the paper.}
We start in section \ref{marginal-distributions} by looking at the limiting marginal distributions 
of the random variables Sine$_\beta[\lambda,\lambda']$ for $\lambda< \lambda'$. 
We first study a classical problem on the exit time of a diffusion trapped in the well of a stationary potential
(obtained by neglecting the slow evolution with time).  
Then, we prove that the jump process of $\lfloor \alpha_\lambda/(2\pi)\rfloor$ converges 
weakly to an (inhomogeneous) Poisson point process by first approximating $\alpha_\lambda$ with 
diffusions processes with piecewise constant drifts on a subdivision of small intervals 
and then by using the convergence of the exit time of the stationary well established previously. 
This requires estimates on the sample paths of a single diffusion, in a spirit similar to \cite{laure-romain,laure}. 
In Section \ref{spatial-independence}, we investigate the asymptotic independence of the numbers of points of Sine$_\beta$ in two disjoint 
intervals.  We prove a crucial estimate regarding the typical relative positions of two diffusions $\alpha_\lambda$ and $\alpha_\lambda'$
for $\lambda < \lambda'$.  Loosely speaking, the main point is to use this estimate to prove that, in the limit $\beta\to 0$, 
the jumps of the process  $\lfloor \alpha_{\lambda}/(2\pi)\rfloor$ immediately follow those of $\lfloor \alpha_{\lambda'}/(2\pi)\rfloor$ 
(see Fig. \ref{without-diff}) while
the processes  $\lfloor \alpha_{\lambda}/(2\pi)\rfloor$ and $\lfloor (\alpha_{\lambda'}-\alpha_{\lambda} )/(2\pi)\rfloor$ never jump at the same time (see Fig. \ref{with-diff}). 
The asymptotic independence follows essentially from the fact that two Poisson point processes adapted to the same filtration are independent 
if and only if they never jump simultaneously.

We gather in the next paragraph important properties of the family $(\alpha_\lambda)$ already established in [Section 2.2, \cite{virag}] that we will use throughout the paper.

\smallskip

\noindent {\it First properties of the coupling of the diffusions $\alpha_{\lambda}$:}
\begin{enumerate}
\item[(i)] For all $\lambda < \lambda'$, $\alpha_{\lambda'} - \alpha_{\lambda}$ has the same distribution as $\alpha_{\lambda'-\lambda}$\,;\smallskip
\item[(ii)] ``Increasing property'':  $\alpha_{\lambda}(t)$ is increasing in $\lambda$\,;\smallskip
\item[(iii)] $\lfloor \alpha_{\lambda}(t)/(2\pi)\rfloor$ is non-decreasing in $t$\,;\smallskip
\item[(iv)] $\E[\alpha_{\lambda}(t)] = \lambda \frac{\beta}{4}\int_0^t e^{-\beta s/4} ds $\,;\smallskip
\item[(v)] $\alpha_\lambda(\infty):=\lim_{t \to \infty} \alpha_{\lambda}(t)/(2\pi)$ exists and is an integer a.s.
\end{enumerate}

\noindent We will also use the following notation:
\begin{align*}
\{x\}_{2\pi} = x - 2\pi \,\Big\lfloor \frac{x}{2\pi}\Big\rfloor\,.
\end{align*}

{\bf Acknowledgments}
Special thanks are addressed to Chris Janjigian and Benedek Valk\'o.
We have benefited from insightful and precise comments from them. Their detailed feedback has helped us to improve the second version
of this manuscript, especially the proofs of Lemmas \ref{def-poissons} and \ref{inclusion-croissante}. 
We are also grateful to them for pointing out references \cite{killip,nakano,kritchevski} 
and the connections with our work. 

We thank St\'ephane Benoist and Antoine Dahlqvist for useful comments and discussions. 

R. A. received funding from the European Research Council under the European
Union's Seventh Framework Programme (FP7/2007-2013) / ERC grant agreement nr. 258237 and thanks the Statslab in DPMMS, Cambridge for its hospitality at the time this work was finished.
The work of L. D. was supported by the 
Engineering and Physical Sciences
Research Council under grant EP/103372X/1 and L.D. thanks the hospitality of the maths department of TU and the Weierstrass institute in Berlin.

\section{Limiting marginal distributions}\label{marginal-distributions}
We are first interested in the limiting law of the random variable Sine$_\beta[0,\lambda]$ when $\beta\to 0$ for a single 
fixed $\lambda$. 
In this case, we re-write the diffusion $\alpha_\lambda$ in a more convenient way: 
\begin{align}\label{diffusion-sine}
d\alpha_\lambda = \lambda \frac{\beta}{4} e^{-\frac{\beta}{4}t} dt + 2 \sin(\frac{\alpha_\lambda}{2}) dB_t \,,
\end{align}
where $B$ is a real standard Brownian motion (which depends on $\lambda$).  
Let us introduce the change of variable $R_\lambda:=\log(\tan(\alpha_{\lambda}/4))$. 
A straightforward computation (see \cite{virag}) shows that:
\begin{align}\label{eq-R}
dR_\lambda= \frac{1}{2} \left(\lambda \frac{\beta}{4} e^{-\frac{\beta}{4}t} \cosh(R_\lambda) +  \tanh(R_\lambda) \right) dt + dB_t\,, 
\quad R_\lambda(0)=-\infty\,.  
\end{align} 

\subsection{Trapping in the stationary potential}
In this subsection, we study an exit time problem for a Langevin diffusion $S_\lambda$ 
evolving in a stationary potential $V_\beta$ defined for $r\in \R$ as 
\begin{align*}
V_\beta(r):= -\frac{1}{2} \left(\lambda \frac{\beta}{4} \sinh(r) +  \log \cosh(r) \right)\,. 
\end{align*}
This problem is relevant to our study thanks to  the slow variation with time  of the non-stationary potential as $\beta\to 0$ in 
which the diffusion $R_\lambda$ evolves. 
The diffusion $S_\lambda$ satisfies the following stochastic differential equation 
\begin{align}\label{eq-S}
dS_\lambda= - V'_\beta(S_\lambda)\,  dt + dB_t\,, \quad S_\lambda(0)=-\infty\,.  
\end{align}
In the whole paper, the diffusions $R_\lambda$ and $S_\lambda$ defined respectively in \eqref{eq-R} and \eqref{eq-S} 
are coupled, driven by the same Brownian motion $B$. Under this coupling, we have almost surely 
\begin{align*}
R_\lambda(t) \le S_\lambda(t)
\end{align*}
for all $t\le \zeta$ where $\zeta$ is the first explosion (stopping) time  
\begin{align*}
\zeta:= \inf\{t\geq 0: S_\lambda(t)=+\infty\}\,. 
\end{align*}

\begin{figure}[h!btp] 
     \center
     \includegraphics[width=10cm]{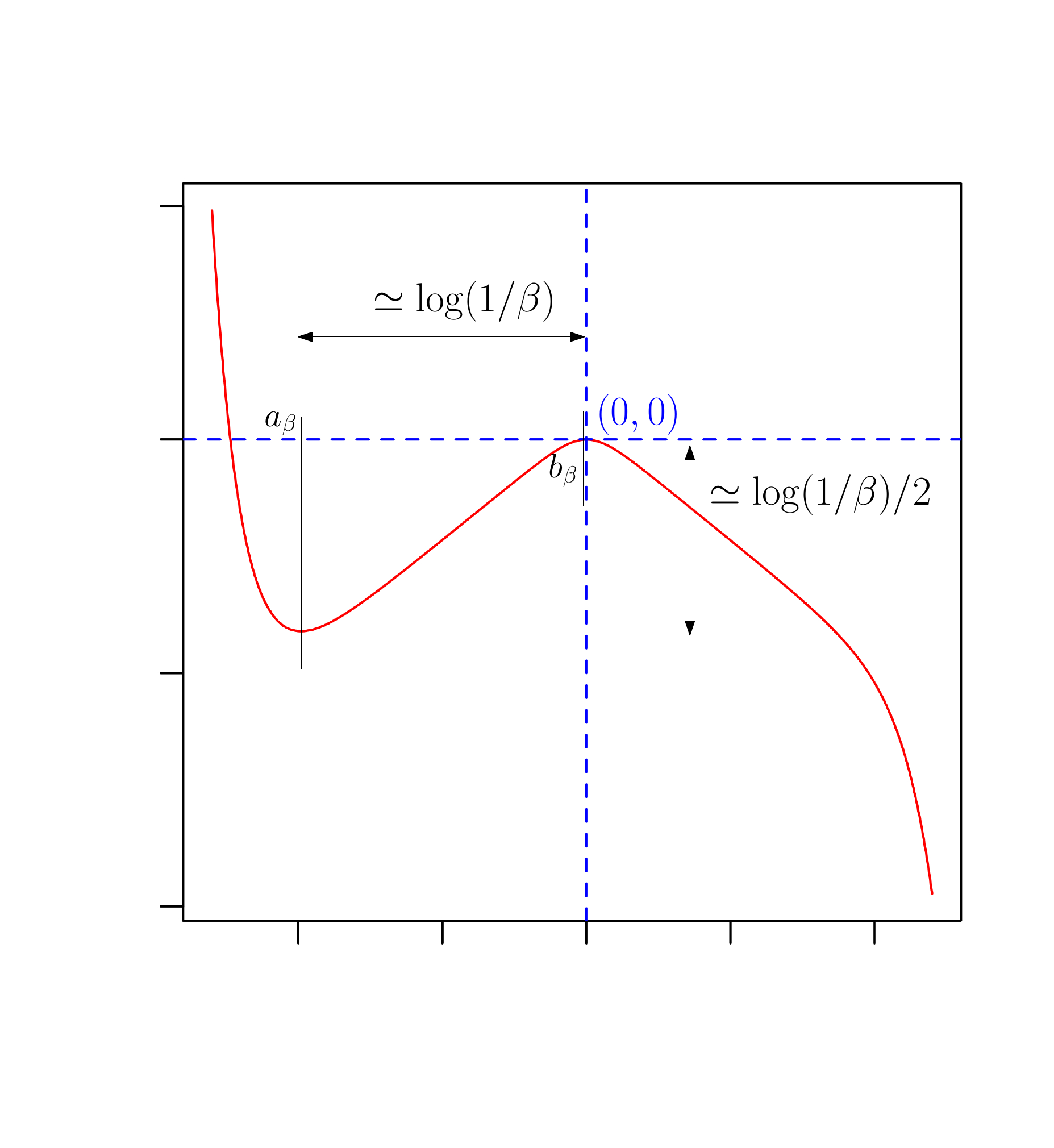}
     \caption{The potential $V_\beta(r)$ as a function of $r$. }\label{fig.potential}
\end{figure}

In this paragraph we investigate the limiting 
law of the stopping time $\zeta$ through its Laplace transform defined for $\xi>0$ as 
\begin{align*}
g_\xi(r) := \E_r[e^{-\xi \zeta}]\,,  
\end{align*}
where $r$ is the initial position of the diffusion 
$S_\lambda$.

We know from classical diffusion theory (see also \cite{laure-romain}) that $g_\xi$ satisfies 
\begin{align}\label{eq-laplace}
\frac{1}{2} g_\xi''(r) - V_\beta'(r) g_\xi'(r) = \xi g_\xi(r)\,, {\mbox{ with }} g_\xi(r)\to 1 {\mbox{ as }}  r\to+\infty\,. 
\end{align}

Let us first examine the expectation of the first explosion time of $S_{\lambda}$. 
Due to the strong barrier separating the local minimum and the local maximum (see Fig. \ref{fig.potential}), 
it is natural to expect the asymptotic of this mean exit time 
 not to depend on the starting point of the diffusion, when $\beta\to 0$, as 
long as it is located in the well (``memory-loss property''). This is the purpose of 
Proposition \ref{mean.time.stat}. We then show that this first exit time properly 
rescaled by its mean value converges in law 
to an exponential distribution (Proposition \ref{conv-exp}) when the starting point is in the well.
\begin{proposition}\label{mean.time.stat}
Suppose that the diffusion $S_\lambda$ starts from $r:= r_\beta$ such that $r_\beta\to -\infty$. Then its expected exit time denoted $t_\beta(r):= \E_r[\zeta]$ has the following equivalent when $\beta \to 0$:
  \begin{align*}
t_\beta(r_\beta) \sim \frac{8\pi}{\beta\lambda}\,.  
\end{align*}
\end{proposition}

{\it Proof of Proposition \ref{mean.time.stat}.}\;
From \eqref{eq-laplace}, it is easy to see that the expected exit time $t_\beta(r)$ satisfies the boundary value problem
\begin{align}\label{eq-expect-exit-time}
\frac{1}{2} t_\beta''(r) - V_\beta'(r) t_\beta'(r) =-1 \,,\quad  {\mbox{ with  }} t_\beta(r)\to 0 {\mbox{ as }}  r\to+\infty\,. 
\end{align}
Solving \eqref{eq-expect-exit-time} explicitly we obtain the integral form 
\begin{align}\label{integral-form-t}
t_\beta(r)= 2 \int_r^{+\infty} dx \int_{-\infty}^x \exp\left(2 \left[V_\beta(x) - V_\beta(y)\right]\right) dy\,. 
\end{align}
By extracting carefully the asymptotic behavior of this integral in the limit $\beta\to 0$ (see Appendix \ref{integrals}), 
we obtained the desired result.  
\qed

The Laplace transform of the exit time satisfies the fixed point equation 
\begin{align*}
g_\xi(r) = 1- 2\xi \int_r^{+\infty} dx \int_{-\infty}^x \exp(2 [V_\beta(x)-V_\beta(y)]) g_\xi(y) dy. 
\end{align*}
With classical arguments similar to those of \cite{laure-romain}, we prove the following proposition in Appendix \ref{integrals}. 
\begin{proposition}\label{conv-exp}
Conditionally on $S_\lambda(0)=r:=r_\beta$ such that $r_\beta \to-\infty$ when $\beta\to 0$, the first explosion time $\zeta$ of the diffusion $(S_\lambda(t))_{t\geq 0}$
converges weakly when rescaled by the  expected exit time $8\pi/(\beta\lambda)$ towards an exponential distribution with mean $1$. 
\end{proposition}

The process $\theta_\lambda$ such that $S_\lambda:=\log\tan(\theta_\lambda/4)$ satisfies
\begin{align*}
d\theta_\lambda= \lambda \frac{\beta}{4} dt + 2 \sin(\frac{\theta_\lambda}{2}) dB_t\,. 
\end{align*}
Proposition \eqref{conv-exp} translates into a convergence of the stopping time $\zeta_{2\pi}:=\inf\{t\geq 0: \theta_\lambda(t)\geq 2\pi\}$. 
 \begin{proposition}\label{conv-exp-theta}
Conditionally on $\theta_\lambda(0)=\theta_0:=\theta_0^\beta$ such that $\theta_0^\beta \to 0_+$ when $\beta\to 0$, the stopping time $
\frac{\beta\lambda}{8\pi}\zeta_{2\pi}$ 
converges weakly  when $\beta\to 0$ towards an exponential distribution with mean $1$. 
\end{proposition}

\subsection{Convergence of the jump process}

We consider the diffusion $\alpha_\lambda$ defined in \eqref{eq-alpha} or equivalently for a single $\lambda$, in \eqref{diffusion-sine}. 
For $k \in \N$, let 
\begin{align}\label{stopping-time}
\zeta_k:=\inf\{t\geq 0: \alpha_\lambda(t) \geq 2 k \pi\}\,. 
\end{align} 
Note that those stopping times correspond to the jumps of the process $\lfloor \alpha_\lambda/(2\pi)\rfloor$.
We denote by $\mathcal{F}_t$ the filtration associated to the diffusion process $\alpha_\lambda$. In this paragraph, we will sometimes omit the subscript $\lambda$ in $\alpha_{\lambda}$ to simplify the notations.
We consider the rescaled empirical measure $\mu_\lambda^\beta$ of the $\frac{\beta}{4} \zeta_k$ defined on $\R_+$  
\begin{align}\label{emp-measure}
\mu_\lambda^\beta[0,t] = \sum_{k=1}^{+\infty} \delta_{\zeta_k}\left[0, \frac{8 \pi t}{\beta}  \right]\,. 
\end{align} 
We divide the time interval $\R_+$ into random small intervals $I_k:=[\frac{S_k}{n}\frac{8\pi}{\beta},\frac{S_{k+1}}{n}\frac{8\pi}{\beta}],k\in \N$, \emph{independent} of the diffusion $\alpha$
where $S_0 = 0$ and:
\begin{align*}
S_k:= \sum_{i=1}^k \tau_i \,, 
\end{align*} 
where the $\tau_i$ are i.i.d. random variables with mean $1$ uniformly distributed on $[1/2,3/2]$.

For each $k\in \N$, conditionally on $\mathcal{F}_{\frac{S_k}{n}\frac{8\pi}{\beta}}$, we define the two diffusion processes $\alpha_{n}^+$ 
and $\alpha_{n}^-$ such that, 
for $t\in I_k$, 
\begin{align*}
d \alpha_{n}^+ = \lambda\frac{\beta}{4} e^{- \frac{S_k}{n}2\pi} dt + 2 \sin(\frac{\alpha_n^+}{2}) dB_t\,, \quad \alpha_{n}^+(\frac{S_k}{n}\frac{8\pi}{\beta}) =  \alpha(\frac{S_k}{n}\frac{8\pi}{\beta})\,,\\
d \alpha_{n}^- = \lambda\frac{\beta}{4} e^{- \frac{S_{k+1}}{n}2\pi} dt + 2 \sin(\frac{\alpha_n^-}{2}) dB_t\,, \quad \alpha_{n}^-(\frac{S_k}{n}\frac{8\pi}{\beta}) =  \alpha(\frac{S_k}{n}\frac{8\pi}{\beta})\,. 
\end{align*}
By the increasing property, it follows that for any $k \in \N$ and $t\in I_k$, we have 
\begin{align*}
 \alpha_{n}^-(t) \le \alpha(t) \le  \alpha_{n}^+(t)\,. 
\end{align*}

\begin{theorem}\label{stopping-times-pp}
As $\beta\to 0$, the empirical measure $\mu_\lambda^\beta$ converges weakly in the space of Radon measure 
(equipped with the topology of vague convergence  \cite{kallenberg}) to a Poisson point process on $\R_+$ 
with inhomogeneous intensity $\lambda \, e^{-2\pi t}\,  dt$. 
In particular, we have, for any $k\in \N$, 
\begin{align*}
\P[\mu_\lambda^\beta[0,t]=k] \to_{\beta\to 0} \exp(-\frac{\lambda}{2\pi}( 1-e^{-2\pi t} ) ) \frac{(\frac{ \lambda}{2\pi})^k  (1-e^{-2\pi t})^k }{k!}\,. 
\end{align*}
\end{theorem}

The proof of Theorem \ref{stopping-times-pp} is done in the next paragraph \ref{estimates-one-diff}. It uses a careful analysis of the behaviour of a single diffusion $\alpha_{\lambda}$ in the small $\beta$ limit.

 Thanks to the equality in law 
\begin{align*}
\mu_\lambda^\beta(\R_+) \stackrel{(d)}{=} \frac{\alpha_\lambda(\infty)}{2\pi}   \stackrel{(d)}{=}  {\rm Sine}_\beta[0,\lambda]\,,
\end{align*}
we easily deduce from Theorem \ref{stopping-times-pp} the convergence of the marginals of the Sine$_\beta$ point process. 
\begin{corollary}\label{marginal-conv-main}
Let $\lambda < \lambda'$. 
The random variable ${\rm Sine}_\beta[\lambda,\lambda']$ converges weakly as $\beta\to 0$ to a Poisson law 
with parameter $\frac{\lambda'-\lambda}{2\pi}$. 
\end{corollary}

\subsection{Estimates for a single diffusion and proof of Theorem \ref{stopping-times-pp}} \label{estimates-one-diff}

We analyse in this paragraph the diffusion $\alpha_{\lambda}$ and the jumps of
$\lfloor \alpha_\lambda/(2 \pi) \rfloor$ in the limit $\beta \to 0$. The results we obtain are derived using the diffusion $R_{\lambda} := \log(\tan(\alpha_{\lambda}/4))$ satisfying the SDE \eqref{eq-R}.
We defer the proof of those technical lemmas in Appendix \ref{proof-auxiliary-results-1}.

\medskip

Lemma \ref{fast-explosion-diff} first states that when the diffusion $\alpha_{\lambda}$ starts just below $2\pi$ modulo $2\pi$, $\lfloor \alpha_\lambda/(2 \pi) \rfloor$  will jump in a short time with probability going to $1$.
\begin{lemma}\label{fast-explosion-diff}
Let $0 <\varepsilon<1$. Conditionally on $\{\alpha_\lambda(0)\}_{2\pi} = 2 \pi - 4 \arctan(\beta^{\varepsilon})$,  
we define the first reaching time 
\begin{equation}\label{reaching-time-mod-2pi}
\zeta^\lambda _{2 \pi  } := \inf\{t\ge 0: \alpha_\lambda(t) = 2 \pi + \lfloor \alpha_{\lambda}(0)/(2\pi) \rfloor \}\,.
\end{equation}
Then, there exists a constant $c >0$ such that for all $\beta>0$ small enough,  
\begin{align*}
\P  \left[ \zeta ^\lambda_{2 \pi } < 9    \log \frac{1}{\beta}  \right]  \geq 1- \beta^c  \,. 
\end{align*}
\end{lemma}

On the large scale-time of the order $1/\beta$, we prove that the time spent by $\alpha_{\lambda}$ near $0$ modulo $2\pi$ for $\alpha_{\lambda}$ is negligible. This is the content of Lemma \ref{expected-time-bad-region-alpha}.

\begin{lemma}\label{expected-time-bad-region-alpha}
Let $t>0$ and 
\begin{align}\label{def-Tbeta}
\Xi_\beta(t) := \E\left[  \int_{0}^{\frac{8\pi }{\beta}t} 1_{\{\{ \alpha_\lambda (u)\}_{2\pi} \ge 4 \arctan(\beta^{1/4}) \}} du \right]  \,. 
\end{align}
Then,  there exists $C>0$ independent of $\beta$ such that 
\begin{align*}
\Xi_\beta(t)\leq \frac{C}{\sqrt{\beta}} \log \frac{1}{\beta} \,. 
\end{align*}
\end{lemma}

\noindent We can now prove Theorem \ref{stopping-times-pp} using the previous estimates:

\medskip

\noindent 
{\it Proof of Theorem \ref{stopping-times-pp}.}

The proof follows the same lines as the proof of [Theorem 4.1 in \cite{laure-romain}]. The idea is to approximate the number of jumps of the diffusion $\alpha$ by those of stationary diffusions and to use the increasing property. To this end, we will use the subdivision of $\R_+$ introduced above and the diffusions $\alpha_n^+$ and $\alpha_n^-$.

From Kallenberg's theorem \cite{kallenberg}, we just need to see that, for any finite 
union $I$ of disjoint and bounded intervals, we have when $\beta\to 0$,  
\begin{align}
\E[\mu_\lambda^\beta(I)] &\longrightarrow \lambda \int_I e^{-2\pi t}\, dt \,,  \label{cv1bis} \\
\P[\mu_\lambda^\beta(I) = 0] &\longrightarrow \exp(- \lambda\int_I e^{-2\pi t}\, dt ) \label{cv2bis} \,. 
\end{align}
Denote by $[t_1;t_2]$ the right most interval of $I$,  
by $J$ the union of disjoint and bounded intervals such that $I=J\cup [t_1;t_2]$ and by $t_0$ the supremum of $J$.

Let us use the random subdivision introduced above. It is crucial to control the position of the diffusion at the starting point of the sub-intervals. As a consequence of Lemma \ref{expected-time-bad-region-alpha}, with large probability, the diffusion $\alpha$ is close to $0$ modulo $2\pi$ in the beginning of each of the sub-intervals overlapping $[0,t_2]$.

More precisely, denote by $\mathcal{C}_k := \{\{\alpha( (8\pi/\beta)S_k/n )\}_{2\pi} \leq 4 \,\arctan(\beta^{1/4})\}$ and consider:
\begin{align*}
\mathcal{C} := \bigcap_{k=0}^{\lfloor 2n t_2\rfloor +1} \mathcal{C}_k.
\end{align*} 
Every sub-interval intersecting $[0,t_2]$ is taken care of in this event as by definition $S_k \in [k/2, 3 k/2]$. Its probability is bounded from above by
\begin{align*}
\P[\mathcal{C}^c] &\leq \sum_{k=0}^{\lfloor 2n t_2\rfloor +1}\P[ \{\alpha((8\pi/\beta) S_k/n)\}_{2\pi} > 4 \,\arctan(\beta^{1/4})]\\
&\leq 2\, \E[\int_{0}^{3n t_2+3} 1_{\{\{\alpha((8\pi/\beta)u/n)\}_{2\pi} > 4 \,\arctan(\beta^{1/4})\}} du]\\
&\leq \frac{n}{4\pi} \beta \; T_{\beta}(0,3t_2+1) \to_{\beta \to 0} 0
\end{align*}
where we used the result of Lemma \ref{expected-time-bad-region-alpha} to have the convergence towards $0$ in the last line.

We now turn to the proof of \ref{cv1bis}. Note that thanks to the linearity of the expectation, we simply need to prove \eqref{cv1bis} for intervals $I$ of the form $I=[0,t_2]$. 
The upper bound simply follows from the SDE form:
\begin{align*}
\E[\alpha(t)] = \lambda \frac{\beta}{4} \int_0^t e^{-\beta s/4} ds \,.
\end{align*}
We immediately derive:
\begin{align*}
\E[\mu_{\beta}[0,t_2]] = \E[\lfloor \alpha((8 \pi/\beta)t_2)/(2\pi) \rfloor] \leq  \lambda \int_0^{t_2} e^{- 2 \pi s} ds \,.
\end{align*}

For the lower bound, denote by $N_k^+$, $N_k^-$ and $N_k$ the number of jumps of $\alpha_n^+$, $\alpha_n^-$ and $\alpha$ in the sub-interval $[(8\pi/\beta)S_k/n, (8\pi/\beta)S_{k+1}/n]$. 
\begin{align*}
\E[\mu_{\beta}[0, t_2]] &\geq \sum_{k=0}^{\lfloor 2n t_2\rfloor +1}\E[N_k^-1_{\{(S_k/n) < t_2\}}| \mathcal{C}_k] - \sum_{k=0}^{\lfloor 2n t_2\rfloor +1} \E[N_k^-1_{\{(S_k/n) < t_2\}}| \mathcal{C}_k] \P[\mathcal{C}_k^c]\,.
\end{align*}
The second sum of the RHS can be simply bounded from above by:
\begin{align*}
\E[\mbox{number of jumps of } \theta_{\lambda} \mbox{ in } [0, 3\, t_2 (8\pi)/\beta]] \times \sup_{k \in \{0, \cdots, \lfloor 2n t_2\rfloor +1\}} \frac{\P[\mathcal{C}^c_k]}{\P[\mathcal{C}_k]}\,.
\end{align*}
The first expectation is bounded from above by a number independent of $\beta$ and the second term (bounded by $\P[\mathcal{C}^c]/(1-\P[\mathcal{C}^c])$) tends to $0$ as $\beta \to 0$.
Moreover, thanks to Proposition \ref{conv-exp-theta}, we have 
\begin{align*}
\E\big[\,N_k^+\, \big{|} \mathcal{C}_k\big] \to_{\beta \to 0} \E[\lambda (\tau_{k+1}/n) \;\exp(-2 \pi S_k/n) ].
\end{align*}
Using the convergence of the Riemann sum when $n \to \infty$, it gives the desired lower bound.

Let us now examine the convergence \ref{cv2bis} for a single interval $[t_1,t_2]$ first. By the Markov property:
\begin{align*}
\P[\mu_{\beta}[t_1,t_2] = 0] \leq \E\Big[\prod_{k=0}^{\infty} \P\big[N_k^+ = 0 \,\big| \mathcal{C}_k, (\tau_i)_i\big] 1_{\{S_{k+1}/n \geq t_1,\;S_{k}/n \leq t_2\}}\Big] + \P[\mathcal{C}^c]\,.
\end{align*}
Using the convergence in each of the sub-intervals (given by Proposition \ref{conv-exp-theta})
\begin{align*}
\P\big[N_k^+ = 0 \, \big{|} \mathcal{C}_k \big] \to_{\beta \to 0} \E[\exp(- \lambda (\tau_{k+1}/n) \;\exp(-2 \pi S_k/n) ))],
\end{align*}
we obtain:
\begin{align*}
\limsup_{\beta \to 0}\P[\mu_{\beta}[t_1,t_2] = 0] \leq \E[\prod_{k=0}^{\infty}\exp(- \lambda (\tau_{k+1}/n) \;\exp(-2 \pi S_k/n) )) 1_{\{S_{k+1}/n \geq t_1,\;S_{k}/n \leq t_2\}}]\,.
\end{align*}
Thanks to the convergence of the Riemann sum when $n \to 0$, we deduce the upper bound. The lower bound can be done using similar techniques than above.
To generalize the result to $I$ finite union of interval, note that thanks to the simple Markov property, we have
\begin{align*}
\P[\mu_\lambda^\beta(I) = 0] = \E[1_{\{\mu_{\beta}(J) = 0\}}\;\P[\mu_{\beta}([t_1,t_2]) = 0| \alpha(t_0)]].
\end{align*}
Iterating the previous argument leads to the result.
\qed

\section{Asymptotic spatial independence}\label{spatial-independence}
\subsection{Ordering of two diffusions}
For $\lambda < \lambda'$, we now control the expected time spent by the process $\{\alpha_{\lambda'}\}_{2 \pi}$ below the process 
$\{\alpha_{\lambda}\}_{2 \pi}$. 
The following Lemma is a crucial step towards the proof of the asymptotic independence of the limiting point process on disjoint intervals and will be used for the proof of Lemmas \ref{inclusion-croissante} and \ref{poissons-indep}. 
 
\begin{lemma}\label{ordering-two-diffusions}
Let $t>0, \lambda < \lambda'$ and 
\begin{align*}
\Theta_\beta(t) :=\E\left[  \int_{0}^{\frac{8\pi }{\beta}t} 1_{\{\{\alpha_{\lambda'}(u)\}_{2\pi} \le\{\alpha_\lambda(u)\}_{2\pi} \}} du\right] \,. 
\end{align*}
Then, there exists two constants $c,C>0$ independent of $\beta$ such that 
\begin{align*}
\Theta_\beta(t)  \leq \frac{C}{\beta^{1-c}} \,. 
\end{align*}

\end{lemma}

{\it Proof.}
We set 
\begin{align*}
\mathcal{E}_u:=\left\{\{\alpha_{\lambda'}(u)\}_{2\pi} <  \{\alpha_\lambda(u)\}_{2\pi}\right\}\,.  
\end{align*}
Before evaluating the probability of the event ${\mathcal{E}_u}$, we need to introduce for $u>0$ 
the last jump of the process $\lfloor\frac{ \alpha_{\lambda'}}{2\pi} \rfloor$ before time $u$, i.e. 
\begin{align*}
\zeta^{(u)}:= \sup \{ k : \zeta_k^{\lambda'} \le u \}\,. 
\end{align*}
The main idea is to prove that any $\mathcal{E}_u$ is associated to a jump time of $\alpha_{\lambda'}$ right before $u$.

\begin{figure}
\centering
\includegraphics[width=9cm]{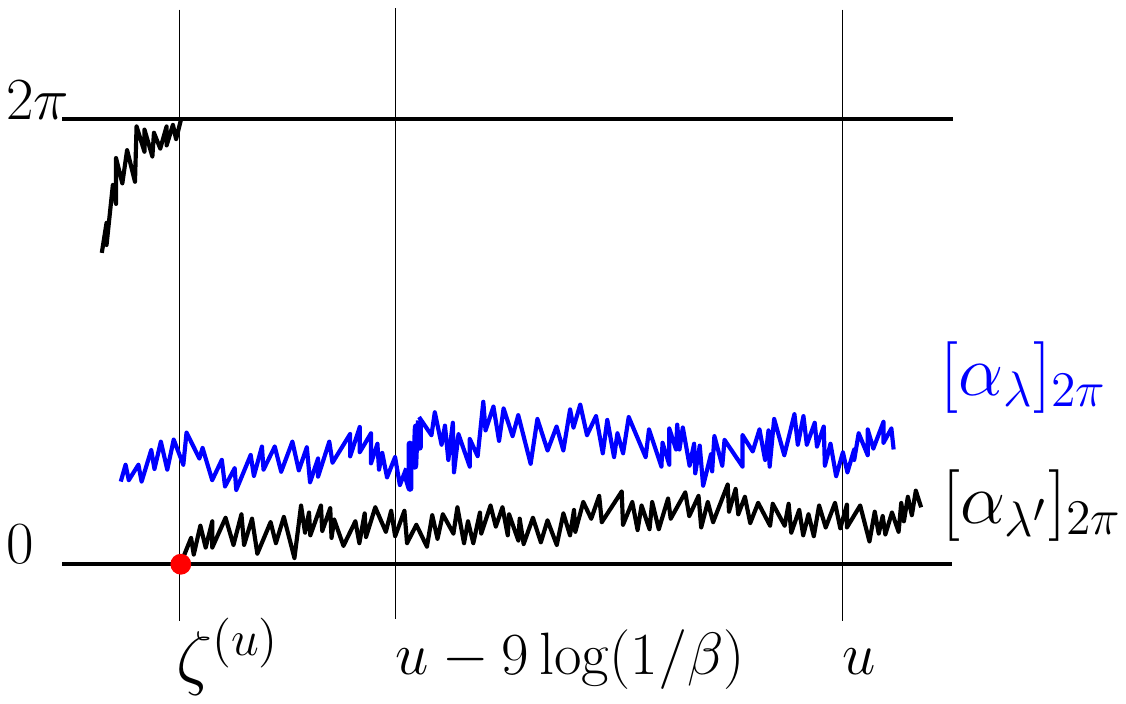}
\caption{Representation of the event $\mathcal{E}_u \cap\{ \zeta^{(u)}  <   u_0 \}$.}\label{fig-alpha-Eu}
\end{figure}

We set $u_0:= u -9 \log\frac{1}{\beta}$ and bound from above the probability of the event ${\mathcal{E}_u}$ by
\begin{align}
\P\left[ \mathcal{E}_u \right]  & \le 
\P\left[ \mathcal{E}_u \cap \left\{\zeta^{(u)} \in \left[u_0, u\right] \right\}  \right] 
+ \P\left[ \mathcal{E}_u \cap\left\{ \zeta^{(u)}  <   u_0  \right\}\right]\notag \\ 
&\le \P\left[ \mathcal{E}_u \cap \left\{\zeta^{(u)} \in \left[u_0, u\right] \right\}  \right] 
+ \P\left[\bigcap_{s\in [u_0,u]} \mathcal{E}_s  \right]  \label{good-shape1}
\end{align}
where we have noticed the inclusion
$\mathcal{E}_u \cap\{ \zeta^{(u)}  <   u_0 \}\subseteq \bigcap_{s\in [u_0,u]} \mathcal{E}_s $ for the second line (see Fig. \ref{fig-alpha-Eu}). Indeed the definition of $\zeta^{(u)} $
implies that the process $\lfloor\frac{ \alpha_{\lambda'}}{2\pi} \rfloor$ has not jumped during the time interval 
$[\zeta^{(u)},u]$ so that the relative ordering (modulo $2\pi$) 
 $\{\alpha_{\lambda'}(t)\}_{2\pi} <  \{\alpha_{\lambda}(t)\}_{2\pi}$ at time $t=u$ has to be preserved for all $t\in [\zeta^{(u)},u)$.

We now tackle the second probability of \eqref{good-shape1} using the fact that $\alpha_\lambda$ is close to $0$ modulo $2\pi$ with large probability:
\begin{align*}
 \P\left[ \bigcap_{s\in [u_0,u]} \mathcal{E}_s  \right] & \le 
  \P\left[ \bigcap_{s\in [u_0,u]} \mathcal{E}_s  \cap 
  \left\{ \{\alpha_\lambda(u_0)\}_{2\pi} \le 4 \arctan(\beta^{1/4}) \right\}  \right]\\ 
  &+   \P\left[\{\alpha_\lambda(u_0)\}_{2\pi} \ge 4 \arctan(\beta^{1/4})   \right]\,. 
\end{align*} 
We introduce the stopping time
\begin{align*}
\zeta_>(u_0) &:= \inf\{t \ge u_0: \{\alpha_{\lambda'}(t)\}_{2\pi} >  \{\alpha_\lambda(t)\}_{2\pi} \} 
\end{align*}
and notice that, conditionally on the event $\mathcal{A}_{u_0} := \mathcal{E}_{u_0} \cap 
  \left\{ \{\alpha_\lambda(u_0)\}_{2\pi} \le 4 \arctan(\beta^{1/4}) \right\}$, we have almost surely for all $t\in [u_0,  \zeta_>(u_0)]$, 
\begin{align*}
\{\alpha_{\lambda'}(t) - \alpha_{\lambda}(t)\}_{2 \pi} &=  2\pi- \left(\{\alpha_{\lambda}(t) \}_{2\pi} - \{\alpha_{\lambda'}(t) \}_{2\pi}\right) \,.  
\end{align*}
Conditionally on the event  $\mathcal{A}_{u_0}$, 
we define a diffusion process $\tilde{\alpha}_{\lambda'-\lambda}$ and its associated first reaching time of $2\pi$ such that
\begin{align*}
d\tilde{\alpha}_{\lambda'-\lambda} (t)&= (\lambda'-\lambda) \frac{\beta}{4}e^{-\frac{\beta}{4}t} dt + {\rm Re}((e^{-i \tilde{\alpha}_{\lambda' -\lambda}(t)} -1)dZ_t)\,,\\  & \mbox{ with }  \tilde{\alpha}_{\lambda'-\lambda}(0)=
2\pi -  ( \{\alpha_\lambda(u_0)\}_{2\pi} -   \{\alpha_{\lambda'}(u_0)\}_{2\pi} )\,, \\ 
& {\mbox{ and }} \tilde{\zeta}_{2\pi} := \inf\{t\ge 0: \tilde{\alpha}_{\lambda'-\lambda} (t) = 2\pi \}\,. 
\end{align*}
Endowed with those definitions, we have the following upper-bound 
\begin{align*}
 \P\left[ \bigcap_{s\in [u_0,u]} \mathcal{E}_s  \cap \mathcal{A}_{u_0}  \right] 
 & = \P\left[\left\{ \zeta_>(u_0) \ge u \right\}  \cap \mathcal{A}_{u_0}  \right]\\ 
 &\le \P\left[\tilde \zeta_{2\pi} \ge 9 \log\frac{1}{\beta} \right] \le \beta^c  
\end{align*}
where we have used the equality in law 
$\tilde{\alpha}_{\lambda'-\lambda}(\cdot) = \alpha_{\lambda'}(u_0+\cdot) - \alpha_\lambda(u_0+\cdot)$  
for the second line 
and Lemma \ref{fast-explosion-diff} (as well as the increasing property) to obtain the last upper-bound. 
Gathering the above inequalities, we obtain 
\begin{align*}
\P\left[ \mathcal{E}_u \right]  & \le
 \P\left[\zeta^{(u)} \in \left[u_0, u\right] \right]  + \beta^{c}  +  \P\left[\{\alpha_\lambda(u_0) \}_{2\pi} \ge 4 \arctan(\beta^{1/4})   \right]\,.
\end{align*}
Now, to conclude, we just have to integrate this latter inequality with respect to $u$. 
First notice that we have almost surely 
\begin{align*}
\int_{9 \log(1/\beta)}^{\frac{8\pi}{\beta}t} 1_{\left\{ \zeta^{(u)} \in \left[u- 9\log\frac{1}{\beta}, u\right] \right \}} du \le 9 \log\frac{1}{\beta} \, 
\mu_{\lambda'}^\beta[0,\frac{8\pi}{\beta}t]\,. 
\end{align*} 
Integrating the other terms as well with respect to $u\in [0,\frac{8\pi}{\beta}t]$ and taking the expectation, we get 
\begin{align*}
\Theta_\beta(t) &\le  9 \log \frac{1}{\beta}\left(1 +\E\left[\mu_{\lambda'}^\beta[0,\frac{8\pi}{\beta}t] \right]\right)  +\Xi_\beta(t) \\
& \le    9 \log \frac{1}{\beta}\left(1 + \frac{\lambda}{2\pi} \right) + 8\pi\,  t\,  \beta^{c-1}  +\Xi_\beta(t) 
\end{align*}
where $\Xi_\beta(t)$ is defined in \eqref{def-Tbeta}.
The conclusion now follows from Lemma \eqref{expected-time-bad-region-alpha}.  
\qed

\subsection{Limiting coupled Poisson processes}
Lemma \ref{ordering-two-diffusions} shall be an important tool 
to prove the asymptotic independence between $\alpha_\lambda(\infty)$ and $\alpha_{\lambda'}(\infty)-\alpha_\lambda(\infty)$ 
for $\lambda < \lambda'$.

Theorem \ref{stopping-times-pp} gives the weak convergence of the random measures  
$\mu_\lambda^\beta$ and $\mu_{\lambda'}^\beta$ in the space of measures on $\R_+$ equipped with the topology of vague convergence denoted $\mathcal{M}_+(\R_+)$.
Due to the equality in law 
\begin{align}\label{equality-in-law}
\alpha_{\lambda'} - \alpha_\lambda \stackrel{(d)}{=} \alpha_{\lambda'-\lambda}\,,  
\end{align}
Theorem \ref{stopping-times-pp} also implies the weak convergence of the (positive) random measure  $\mu_{\lambda'-\lambda}^\beta$
such that for all $t\ge 0$, 
\begin{align*}
\mu_{\lambda'-\lambda}^\beta[0,t] := \lfloor \frac{\alpha_{\lambda'}(t)-\alpha_\lambda(t)}{2\pi} \rfloor
\end{align*}
towards a Poisson measure $\mathcal{P}_{\lambda'-\lambda}$ with intensity $(\lambda'-\lambda)\,e^{-2\pi t}dt$. 

We now work with the two diffusions $\alpha_\lambda$ and $\alpha_{\lambda'}$ for $\lambda < \lambda'$ \emph{coupled 
according to \eqref{eq-alpha}}. 
We are interested in the limiting joint distribution of the triplet of random measures $(\mu_{\lambda}^\beta, \mu_{\lambda'}^\beta, \mu_{\lambda'-\lambda}^\beta)$ according to this coupling.

From the above convergences, it is straightforward to deduce the relative-compactness of the family of the triplets of (random) measures
\begin{equation}\label{family-meas-beta}
\{ \big(\mu_\lambda^\beta, \mu_{\lambda'}^\beta,  \mu^\beta_{\lambda''-\lambda'}\big), \beta>0  \} 
\end{equation}
for the weak topology over $(\mathcal{M}_+(\R_+))^3$ equipped with the product topology of vague convergence.
 
Let us take a sequence $\beta_k \to 0$ when $k\to\infty$ such that the processes
\begin{equation}\label{joint-conv}
\big( \mu_\lambda^{\beta_k}, \mu_{\lambda'}^{\beta_k},\mu_{\lambda''-\lambda'}^{\beta_k} \big)
\end{equation} 
converge jointly weakly in the product space when $k\to \infty$ to a triplet 
\begin{align*}
(\mathcal{P}^{(s)}_{\lambda}, \mathcal{P}^{(s)}_{\lambda'}, \mathcal{P}^{(s)}_{\lambda''-\lambda'})
\end{align*}
of point measures on $\R_+$ whose marginal distributions are given respectively by 
the law of the Poisson measures $\mathcal{P}_{\lambda}$ and $\mathcal{P}_{\lambda'}$ 
and $\mathcal{P}_{\lambda''-\lambda'}$ and whose joint law depends \emph{a priori} on the chosen sub-sequence $(s):=(\beta_k)$. In the following, we will study this triplet and we will drop the superscript $(s)$ to ease the notations. 
We shall in fact see later that all the possible limit point have the same law. 
Therefore the law of the triplet does not depend on $(s)$ and the weak convergence of \eqref{family-meas-beta} holds (see Remark \ref{remark-inclusion-facult}). In the next Lemma, we regard the point measures $\mathcal{P}_{\lambda},\mathcal{P}_{\lambda'}$ 
and $\mathcal{P}_{\lambda''-\lambda'}$ as Poisson point processes and prove that they are indeed jointly Poisson processes 
on a common filtration.  This is an important step for our needs.

\begin{lemma}\label{def-poissons}
Let $\lambda < \lambda' < \lambda''$
and $\mathcal{F}:=(\mathcal{F}_t)_{t\ge 0}$ be the natural filtration associated to the process 
$(\mathcal{P}_{\lambda}, \mathcal{P}_{\lambda'}, \mathcal{P}_{\lambda''-\lambda'})$ i.e. such that 
\begin{equation*} \mathcal{F}_t:=\sigma ( \mathcal{P}_{\lambda}(s), \mathcal{P}_{\lambda'}(s), \mathcal{P}_{\lambda''-\lambda'}(s) ,0 \le s\le t).\end{equation*}
Then, the processes  $\mathcal{P}_{\lambda}$ and $\mathcal{P}_{\lambda'}$ and $\mathcal{P}_{\lambda''-\lambda'}$ are 
$(\mathcal{F}_t)$-Poisson point processes with respective intensities 
$\lambda e^{-2\pi t}dt$ and $\lambda' e^{-2\pi t}dt $ and $(\lambda''-\lambda') e^{-2\pi t}dt$.
\end{lemma}  
The main point used to prove this Lemma is that all the diffusions are measurable 
with respect to the same driving (complex) Brownian motion $(Z_t)$ 
as they are strong solutions of the SDEs \eqref{eq-alpha}.   
A proof can be found in Appendix \ref{appendix-def-poissons}.

To fix notations, we will denote by $(\xi_i^\lambda)_{i \in \N}$ resp. $(\xi_i^{\lambda'})_{i \in \N},(\xi_i^{\lambda'-\lambda})_{i \in \N}$ the 
points associated to the Poisson processes
 $\mathcal{P}_{\lambda}$ and $\mathcal{P}_{\lambda'}$ and $\mathcal{P}_{\lambda'-\lambda}$ such that 
 \begin{align*}
\mathcal{P}_{\lambda}[0,t] = \sum_i \delta_{\xi_i^\lambda}[0,t]\,, \quad \mathcal{P}_{\lambda'}[0,t] = \sum_i \delta_{\xi_i^{\lambda'}}[0,t]
\,, \quad \mathcal{P}_{\lambda'-\lambda}[0,t] = \sum_i \delta_{\xi_i^{\lambda'-\lambda}}[0,t] \,. 
\end{align*}
For $\beta >0$, we also recall the notations
\begin{align*}
\mu^\beta_{\lambda}[0,t] = \sum_i \delta_{\zeta_i^\lambda}[0,\frac{8 \pi t}{\beta}]\,, \quad \mu^\beta_{\lambda'}[0,t] = \sum_i \delta_{\zeta_i^{\lambda'}}[0,\frac{8 \pi t}{\beta}]
\,, \quad \mu^\beta_{\lambda'-\lambda}[0,t] = \sum_i \delta_{\zeta_i^{\lambda'-\lambda}}[0,\frac{8 \pi t}{\beta}] \,.
\end{align*} 

\begin{lemma}\label{inclusion-croissante}
Let $\lambda < \lambda'$. Then, we have almost surely
\begin{align*}
\mathcal{P}_{\lambda} \subseteq \mathcal{P}_{\lambda'}
\end{align*}
 i.e. 
for all $ i \in \N$, there exists $j\ge i$ such that $\xi_i^\lambda = \xi_j^{\lambda'}$. 
\end{lemma}

{\it Proof.}
We have to prove that for any $t>0$, 
\begin{align}\label{proba-not-in}
\P\left[ \exists i: \xi_i^\lambda < t, \xi_i^{\lambda}\not\in \mathcal{P}_{\lambda'} \right] = 0\,. 
\end{align}
The probability \eqref{proba-not-in} is the increasing limit of the sequence $(p_n)_{n\in\N}$ defined as 
\begin{align}\label{proba-not-in-n}
p_n:=\P\left[ \exists i: \xi_i^\lambda < t, \forall j\ge i, | \xi_i^{\lambda} - \xi_j^{\lambda'}|> \frac{1}{2n} \right] \,. 
\end{align}
It suffices to prove that $p_n=0$ for any $n$. 
Let us introduce the probability   
\begin{align}\label{discrete-proba1}
p_n^\beta:= 
\P\left[ \exists i: \frac{\beta}{8 \pi}\zeta_i^\lambda <  t, \forall j\ge i, \frac{\beta}{8 \pi}| \zeta_i^{\lambda} - \zeta_j^{\lambda'}|> \frac{1}{2n} \right]\,.  
\end{align}
Denote by $\mathcal{M}_p(\R_+)$ the space of point measures and recall that it is 
closed in the space $ \mathcal{M}_+(\R_+)$ for the vague convergence topology.
Notice that the set 
\begin{align*}
\{(\mu,\nu) \in (\mathcal{M}_p(\R_+))^2 \;:\; \mu = \sum_{i=1}^{\infty} \delta_{x_i}, \nu =  \sum_{i=1}^{\infty} \delta_{y_i}, \exists i \mbox{ such that } y_i <t \mbox{ and } \forall j, \;|y_i - x_j| > \frac{1}{2n}\}
\end{align*}
is open in $(\mathcal{M}_p(\R_+))^2$ equipped with the product topology of the vague convergence on the space of point measures. It comes from the straightforward fact that if 
 $\mu_k \in \mathcal{M}_p(\R_+)$ converges towards $\mu \in \mathcal{M}_p(\R_+)$ for the vague topology, the points of $\mu_k$ belonging to
 $[0,t]$ for all but finitely many $k$ converge to points of $\mu$ in $[0,t]$. 

If $n$ is fixed,  we therefore have  $p_n \le \liminf_{\beta\to 0} p_n^\beta $ using 
the joint convergence of $(\mu_\lambda^\beta,\mu_{\lambda'}^\beta)$ in the space $(\mathcal{M}_p(\R_+))^2$ 
(along the subsequence $\ell$) and the Portmanteau theorem. 
It suffices to check that
\begin{align*}
  \limsup_{\beta\to 0} p_n^\beta = 0\,. 
\end{align*} 
We need to work with a random subdivision of the interval $[0, \frac{8 \pi t}{\beta}]$. As before, we consider a sequence 
$(\tau_k)_{k \in \N}$
of i.i.d. positive random variables distributed uniformly in $[\frac{1}{2},\frac{3}{2}]$ 
and form the sum $S_k=\sum_{i=1}^k \tau_i$. 

Noting that for all $x,y$ such that $|x-y|\le 1/(2n)$, there exists $k\in \N$ such that $x,y \in[ S_k/n ,S_k/n + 
2/n ]$, we obtain 
\begin{align*}
p_n^\beta \le  \P\left[  \exists  k \le \lfloor 2nt \rfloor +1 : \lfloor \frac{ \alpha_\lambda}{2\pi}  \rfloor \mbox{ jumps on the interval }   
\frac{8 \pi}{\beta} \left[ \frac{S_k}{n} ,  (\frac{S_k}{n} + \frac{2}{n} )\right] \mbox{ but not }  \lfloor \frac{ \alpha_{\lambda'}}{2\pi}  \rfloor \right]\,. 
\end{align*}
Due to the increasing property, the event inside the probability can not happen if the process
$\{\alpha_\lambda\}_{2\pi}$ starts below $\{\alpha_{\lambda'}\}_{2\pi}$ at the beginning of the interval. 
Therefore,
\begin{align*}
p_n^\beta \le   
\sum_{k=1}^{\lfloor 2 n t\rfloor+1} \P\left[ \{\alpha_{\lambda'}( \frac{8 \pi}{\beta} \frac{S_k}{n} ) \}_{2\pi} \le[ \alpha_{\lambda}( \frac{8 \pi}{\beta} \frac{S_k}{n} ) ] _{2\pi}\right] \,. 
\end{align*}
which can in turn be upper-bounded as follows
\begin{align}
\sum_{k=1}^{\lfloor 2 n t\rfloor+1} \P\left[ \{\alpha_{\lambda'}( \frac{8 \pi}{\beta} \frac{S_k}{n} ) \}_{2\pi} \le[ \alpha_{\lambda}( \frac{8 \pi}{\beta} \frac{S_k}{n} ) ] _{2\pi}\right] & \le \frac{\beta n}{4\pi} \,  \E\left[\int_0^{\frac{8\pi}{\beta} (3 t+1 ) } 1_{\left\{\{\alpha_{\lambda'}(  u) \}_{2\pi} \le[ \alpha_{\lambda}( u ) ] _{2\pi}\right\}} du \right] \notag\\ 
&\le \frac{\beta n}{4\pi} \Theta_\beta(3t+1) \, \le \beta^c   \label{estimate1}
\end{align}
where we have used Lemma \ref{ordering-two-diffusions} to obtain the last inequality which holds for $\beta$ small enough
( $c$ is a constant which does not depend on $\beta$).   

\qed

\begin{lemma}\label{poissons-indep}
Let $\lambda < \lambda'$. Then the $(\mathcal{F}_t)$-Poisson point processes  $\mathcal{P}_\lambda$ and $\mathcal{P}_{\lambda'-\lambda} $ are independent.
\end{lemma}

{\it Proof of Lemma \ref{poissons-indep}.}
From a classical result (see Proposition (1.7) Chapter XII, \S 1, p.473 in \cite{revuz-yor}) on Poisson processes, we know that it suffices to prove that 
the two $(\mathcal{F}_t)$-Poisson processes $\mathcal{P}_\lambda$ and $\mathcal{P}_{\lambda'-\lambda} $ 
do not jump simultaneously, i.e. that for $t>0$, 
\begin{align}\label{no-simultaneous-jump}
\P\left[\exists i, j \in \N :  \xi_i^\lambda < t, \, \xi_j^{\lambda'-\lambda} < t, \, \xi_i^\lambda = \xi_j^{\lambda'-\lambda} \right] = 0\,. 
\end{align}
For $n \in \N$, we consider the probability  
\begin{align}\label{discrete-proba}
p_n^\beta:=\P \left[\exists i, j \in \N : \frac{\beta}{8 \pi}\zeta_i^\lambda < t, \, \frac{\beta}{8 \pi}\zeta_j^{\lambda'-\lambda} < t, \, 
\frac{\beta}{8 \pi} | \zeta_i^\lambda - \zeta_j^{\lambda'-\lambda} | < \frac{1}{2n}  \right]\,. 
\end{align}
If $n$ is fixed, then we have the convergence (the studied set is open as in the proof Lemma \ref{inclusion-croissante})
\begin{align*}
\liminf_{\beta\to 0} p_n^\beta \ge p_n:=  
\P\left[\exists i, j \in \N :  \xi_i^\lambda < t, \, \xi_j^{\lambda'-\lambda} < t, \, |\xi_i^\lambda - \xi_j^{\lambda'-\lambda} | <  \frac{1}{2n} \right]\,. 
\end{align*}
To prove \eqref{no-simultaneous-jump}, it therefore suffices to prove that 
\begin{align}\label{double-limit}
\limsup_{n \to \infty} \limsup_{\beta\to 0} p_n^\beta =0\,. 
\end{align}
For this, we need to work with a random subdivision of the interval $[0, \frac{8 \pi t}{\beta}]$. As before, we consider a sequence 
$(\tau_k)_{k \in \N}$
of i.i.d. positive random variables uniformly distributed in $[\frac{1}{2},\frac{3}{2}]$ and independent of the processes
and form the sum $S_k=\sum_{i=1}^k \tau_i$ ($S_0 :=0$).

\begin{figure}
\centering
\includegraphics[width=9cm]{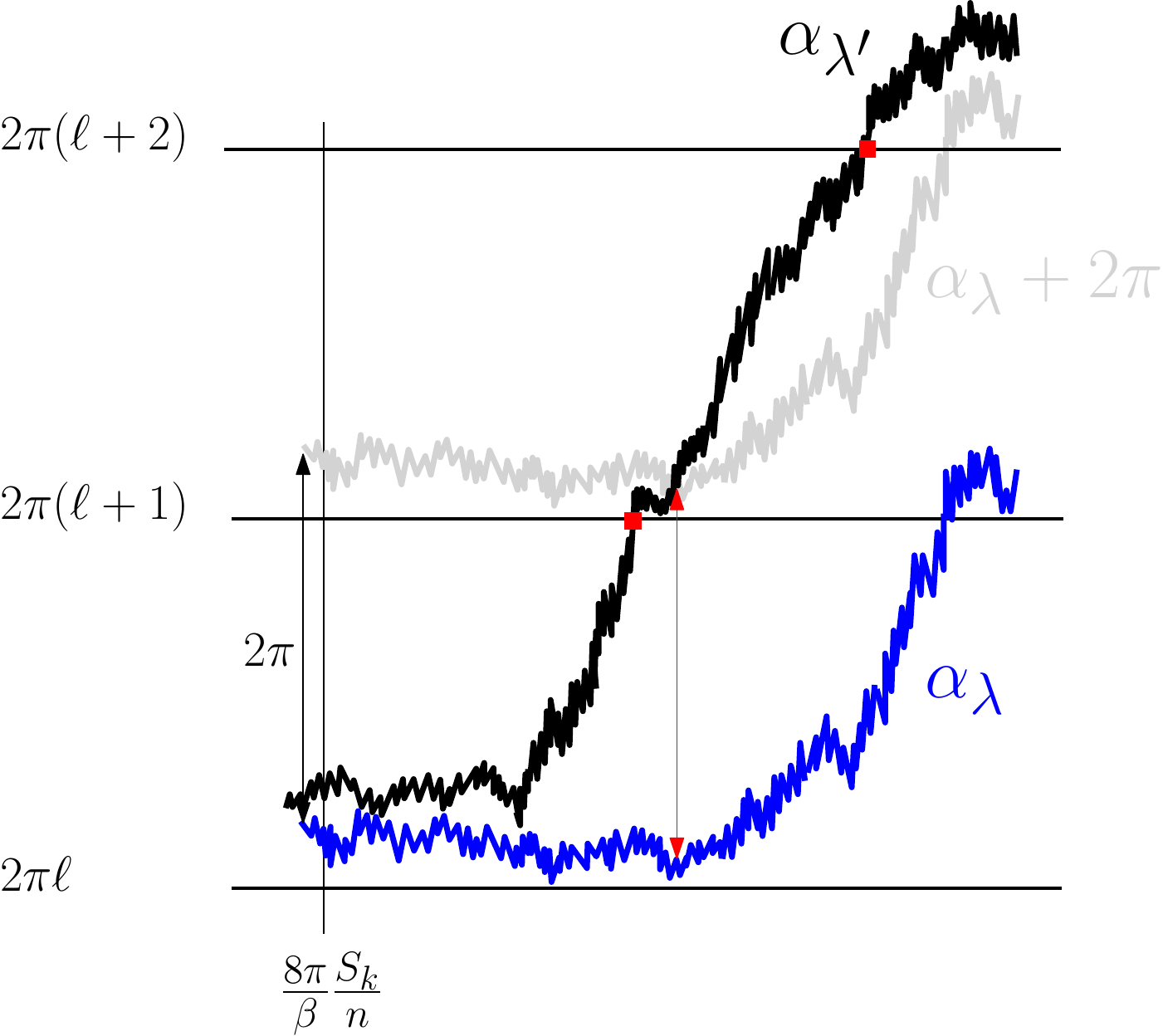}
\caption{Illustration for inclusion \eqref{event1}.}\label{fig-2jumps}
\end{figure}

The probability \eqref{discrete-proba} can be upper-bounded  as follows 
\begin{align}
&\P\left[\exists i,j \in \N, k \le \lfloor 2 n t\rfloor+1  :  \frac{8 \pi}{\beta} \frac{S_k}{n} \le \zeta_i^\lambda, \zeta_j^{\lambda'-\lambda} \le 
\frac{8 \pi}{\beta} (\frac{S_k}{n} +\frac{2}{n} ) ] \right]\notag  \\
&= \P\left[\exists k \le \lfloor 2 n t\rfloor+1 : \lfloor \frac{\alpha_{\lambda'} - \alpha_\lambda}{2\pi}  \rfloor \mbox{ and }  \lfloor \frac{ \alpha_\lambda}{2\pi}  \rfloor \mbox{ both jump on the interval }   \frac{8 \pi}{\beta} \left[ \frac{S_k}{n} ,  (\frac{S_k}{n} + \frac{2}{n} )\right] \right]\notag\\ 
\label{sums-to-estimate}
&\le \sum_{k=1}^{ \lfloor 2 n t\rfloor +1} \P\left[ \{\alpha_{\lambda'}( \frac{8 \pi}{\beta} \frac{S_k}{n} ) \}_{2\pi} \le[ \alpha_{\lambda}( \frac{8 \pi}{\beta} \frac{S_k}{n} ) ] _{2\pi}\right] \\ &+ \sum_{k=1}^{\lfloor 2 n t\rfloor+1 }\P\left[ \lfloor \frac{ \alpha_{\lambda'}}{2\pi}  \rfloor  \mbox{ jumps two times during the interval }  \frac{8 \pi}{\beta} \left[ \frac{S_k}{n} ,  (\frac{S_k}{n} + \frac{2}{n} )\right] \right]\,. 
\end{align}
For this bound, we have used the fact that, conditionally on 
\begin{equation}\label{ordering-condition}
\{\alpha_{\lambda}( \frac{8 \pi}{\beta} \frac{S_k}{n} ) \}_{2\pi} \le[ \alpha_{\lambda'}( \frac{8 \pi}{\beta} \frac{S_k}{n} )] _{2\pi}\,, 
\end{equation}
 the increasing property and the equality in law \eqref{equality-in-law} 
impose 
\begin{align}
&\left\{  
\lfloor \frac{\alpha_{\lambda'} - \alpha_\lambda}{2\pi}  \rfloor  \mbox{ and }  \lfloor \frac{ \alpha_\lambda}{2\pi}  \rfloor 
\mbox{ both jump on the interval }   \frac{8 \pi}{\beta} \left[ \frac{S_k}{n} ,  (\frac{S_k}{n} + \frac{2}{n} )  \right]
\right \} \label{event1} \\
&\subseteq 
\left\{ \lfloor \frac{ \alpha_{\lambda'}}{2\pi}   \rfloor  \mbox{ jumps two times during the interval }  \frac{8 \pi}{\beta} \left[ \frac{S_k}{n} ,  (\frac{S_k}{n} + \frac{2}{n} )\right]  \right \} \notag \,. 
\end{align}
Indeed, the equality in law \eqref{equality-in-law} implies that $\lfloor \frac{\alpha_{\lambda'}(t)-\alpha_\lambda(t)}{2\pi} \rfloor$ is increasing with respect to $t$ (once the difference $\alpha_{\lambda'}-\alpha_\lambda$ 
has reached the value $2 k \pi$ where $k\in \N$, it remains forever above this value).  
This fact and the increasing property imply that under the event \eqref{event1}, 
the process  $ \lfloor \frac{ \alpha_{\lambda'}}{2\pi}   \rfloor$ has to jump two times on the interval $\frac{8 \pi}{\beta} [ \frac{S_k}{n} ,  
(\frac{S_k}{n} + \frac{2}{n}) ]$ (see Fig. \ref{fig-2jumps}). 
We now estimate the two sums in \eqref{sums-to-estimate}. 
The first one was already tackled in \eqref{estimate1}.

For the other sum in \eqref{sums-to-estimate}, we write 
\begin{align*}
\sum_{k=1}^{\lfloor 2 n t\rfloor+1 }& \P\left[ \lfloor \frac{ \alpha_{\lambda'}}{2\pi}  \rfloor  \mbox{ jumps two times during the interval }  \frac{8 \pi}{\beta} \left[ \frac{S_k}{n} ,  (\frac{S_k}{n} + \frac{2}{n} )\right] \right] \\ &= 
\sum_{k=1}^{\lfloor 2 n t\rfloor+1 }\P\left[ \mu_{\lambda'}^\beta\left[\frac{8 \pi}{\beta} \left[ \frac{S_k}{n}  , ( \frac{S_k}{n}   + \frac{2}{n} )\right]\right]  =2 \right]\,. 
\end{align*}
If $n$ is fixed and $\beta\to 0$, Theorem \ref{stopping-times-pp} gives the following convergence 
\begin{align}
\sum_{k=1}^{\lfloor 2 n t\rfloor+1 }&\P\left[ \mu_{\lambda'}^\beta[\frac{8 \pi}{\beta} \frac{S_k}{n} [1 ,  1 + \frac{2}{n} ]] =2 \right]\notag  \\
&\rightarrow_{\beta\to 0} 
\frac{\lambda^2}{16\pi^2}\sum_{k=1}^{\lfloor 2 n t\rfloor+1 } \exp\left(-\frac{\lambda}{2\pi} e^{-2\pi\frac{S_k}{n}} (1-e^{-2\pi\frac{2}{n}}) \right)
e^{-4\pi\frac{S_k}{n}} (1-e^{-2\pi\frac{2}{n}}) ^2 \notag \\ & = O(n^{-1}) \label{estimate2} \,. 
\end{align}
The convergence \eqref{double-limit} now follows from the two estimates \eqref{estimate1} and \eqref{estimate2} and the independence between 
the two Poisson processes $\mathcal{P}_{\lambda}$ and $\mathcal{P}_{\lambda'-\lambda}$ is proved. 
\qed

\begin{rem}\label{remark-inclusion-facult}
Using similar arguments as in the proof of Lemma \ref{inclusion-croissante}, one could also show that for all $\lambda < \lambda'$,
$\mathcal{P}_{\lambda'-\lambda} \subseteq \mathcal{P}_{\lambda'}$. Using then Lemma \ref{poissons-indep}, it implies the almost-sure equality $\mathcal{P}_{\lambda'} = \mathcal{P}_{\lambda} +\mathcal{P}_{\lambda'-\lambda}$\footnote{Note that this equality is obviously wrong when $\beta >0$ i.e we don't have the equality between $\lfloor\alpha_{\lambda'}(t)/(2\pi)\rfloor$ and $\lfloor\alpha_{\lambda'}(t)/(2\pi)\rfloor + \lfloor(\alpha_{\lambda'}(t) - \alpha_{\lambda}(t))/(2\pi)\rfloor$.}. It therefore determines the law of the triplet $(\mathcal{P}_{\lambda},\mathcal{P}_{\lambda'},\mathcal{P}_{\lambda'-\lambda})$ and the weak convergence of the family \eqref{family-meas-beta} when $\beta \to 0$ holds.

Nevertheless, we do not need this stronger result as we only use for the $Sine_{\beta}$ case the previous equality for the total number of points in $\R_+$ i.e. $\alpha_{\lambda'}(\infty)/(2\pi) = \alpha_{\lambda}(\infty)/(2\pi) + (\alpha_{\lambda'}(\infty) - \alpha_{\lambda}(\infty))/(2\pi)$
\end{rem}

As an immediate consequence of Lemma \ref{inclusion-croissante} and Lemma \ref{poissons-indep}, we obtain the following result 
using the characterization of independent Poisson processes mentioned above.  
\begin{lemma}\label{final-independence}
Let $\lambda_0<\lambda_1<\lambda_2$. Then, the $(\mathcal{F}_t)$-Poisson processes
$\mathcal{P}_{\lambda_0}$ and $\mathcal{P}_{\lambda_2-\lambda_1}$ are independent. 
\end{lemma}

\subsection{Proof of Theorem \ref{main}}

From Kallenberg's theorem \cite{kallenberg}, we simply need to check that, for any finite 
union $I$ of disjoint and bounded intervals, we have when $\beta\to 0$,  
\begin{align}
\E[{\rm Sine}_\beta(I)] &\longrightarrow  \frac{1}{2\pi} \int_I  d\lambda \,,  \label{cv1bis-main} \\
\P[{\rm Sine}_\beta(I) = 0] &\longrightarrow \exp(- \frac{1}{2\pi}  \int_I  d\lambda ) \label{cv2bis-main} \,. 
\end{align}
The convergence \eqref{cv1bis-main} follows from Corollary \ref{marginal-conv-main}.

Towards \eqref{cv2bis-main}, we consider first a union of two disjoint interval of the form 
\begin{align*}
I:= [\lambda_0,\lambda_1] \cup [\lambda_2,\lambda_3]\,,
\end{align*}
where $\lambda_0 < \lambda_1 < \lambda_2 < \lambda_3$. 
By translation invariance of the Sine$_\beta$ point process, we can suppose without loss of generality  that $\lambda_0=0$.  
By definition of Sine$_\beta$, we have the following equality in law 
\begin{align*}
\left({\rm Sine}_\beta[0,\lambda_1], {\rm Sine}_\beta[\lambda_2, \lambda_3], \right) \stackrel{(d)}{=} 
\left( \frac{\alpha_{\lambda_1}(\infty)}{2\pi} , \frac{\alpha_{\lambda_3}  (\infty)- \alpha_{\lambda_2} (\infty) }{2\pi} \right) \,. 
\end{align*}
Using Lemma \ref{final-independence}, we deduce that 
\begin{align*}
\P&\left[ {\rm Sine}_\beta[0,\lambda_1] = 0 , {\rm Sine}_\beta[\lambda_2, \lambda_3] =0\right]\\ 
&\rightarrow_{\beta\to 0} \P\left[\mathcal{P}_{\lambda_1}(\infty) = 0 \right] \P\left[\mathcal{P}_{\lambda_3-\lambda_2}(\infty) = 0 \right] 
= \exp\left(- \frac{\lambda_1}{2\pi} \right)\exp\left(- \frac{\lambda_3-\lambda_2}{2\pi} \right) \,. 
\end{align*}

Now, if $k \ge 3$ and $I :=[\lambda_0,\lambda_1] \cup \cdots \cup [\lambda_{k-1},\lambda_k]$ where $\lambda_0 \le \cdots \le \lambda_k$, the limiting Poisson point processes $(\mathcal{P}_{\lambda_i - \lambda_{i-1}})_{i = 1,\cdots,k}$ obtained as the limit of the processes $\lfloor (\alpha_{\lambda_i} - \alpha_{\lambda_{i-1}})/(2\pi)\rfloor$ never jump simultaneously (there are pairwise independent). Therefore, it is a family of independent processes  (see Proposition (1.7) Chapter XII, \S 1, p.473 in \cite{revuz-yor}).

Theorem \ref{main} is proved.
\qed

\appendix 

\section{Proof of auxiliary results.}\label{proof-auxiliary-results}
\subsection{Asymptotic of integrals and Laplace transform}\label{integrals}

{\it Asymptotic of the integral \eqref{integral-form-t}.}
We denote by $a_\beta$ (respectively $b_\beta$) the local minimum (resp. maximum) of the potential $V_\beta(r)$ and derive asymptotic expansions for $a_\beta$ and $b_\beta$ in the limit $\beta\to 0$ using
\begin{align*}
V_\beta'(r)=0 \Leftrightarrow \frac{\lambda\beta}{8} = \frac{e^{-r}-e^{r}}{(e^r+e^{-r})^2}\,. 
\end{align*}
We deduce that 
\begin{align*}
a_\beta&= \log \beta + \log(\lambda) - 3 \log 2 + O(\beta^2)\,,\\
b_\beta&= - \frac{\lambda \beta}{4} + o(\beta^2)\,,\\
V_\beta(a_\beta)&= \frac{1}{2} \log \beta + \frac{1}{2} \log \lambda - \log 2 + \frac{1}{2} + O(\beta^2)\,, \\
V_\beta(b_\beta)&= \frac{\lambda^2}{64} \beta^2 + O(\beta^2) = O(\beta^2)\,. 
\end{align*}

More generally, we have for $x, y \in \R$ 
\begin{align*}
V_\beta(a_\beta+y) &= \frac{1}{2} \log\frac{\beta\lambda}{4} + \frac{y+e^{-y}}{2} - \frac{1}{2} (\frac{\lambda\beta}{8})^2 e^{y}
- \frac{1}{2} \log (1+(\frac{\lambda\beta}{8})^2 e^{2 y} ) +  O(\beta^2)\,, \\
V_\beta(b_\beta+x) & = - \frac{1}{2} \log \frac{e^x+e^{-x}}{2} (1-\frac{\lambda\beta}{8} + O(\beta^2)) - \frac{\lambda\beta}{8} \sinh(x)
+\frac{\lambda^2 \beta^2}{64}\cosh(x) + O(\beta^2)\,.  
\end{align*}

Those computations permit us to find the asymptotic behaviour of $t_\beta(r_\beta)$ 
(using also the bounded convergence 
theorem), 
\begin{align}\label{equiv-t-beta}
t_\beta(r_\beta) \sim \frac{8}{\beta\lambda} \int_{r_\beta-b_\beta}^{+\infty} \frac{dx}{\cosh(x)} \int_{-\infty}^{x+ a_\beta} 
e^{-y} e^{-e^{-y}} dy\sim  \frac{8}{\beta\lambda} \int_{r_\beta-b_\beta}^{+\infty} \frac{dx}{\cosh(x)} \sim \frac{8\pi}{\beta\lambda} \,.  
\end{align}

\qed

{\it Proof of Proposition \ref{conv-exp}.}

We consider the Laplace transform $g_{\beta\lambda/(8 \pi) \xi}(r)$ of the rescaled exit time $\beta\lambda/(8 \pi) \zeta$
of the diffusion starting from $r$. 
If suffices to prove that, if  $(r_\beta)_{\beta>0}$ is such that $r_\beta\to -\infty$ as $\beta\to 0$, we have  for any $0 <\xi <1$, 
$g_{\beta\lambda/(8 \pi) \xi}(r_\beta)\to 1/(1+\xi)$ (Laplace transform of an exponential distribution with parameter $1$).  
To simplify notations, 
we will just write $g_\beta(r)$ for $g_{\beta\lambda/(8 \pi) \xi}(r)$.  

Let us recall the fixed point equation 
\begin{align}\label{fixed-point-rescaled}
g_\beta(r) = 1-2 \frac{\beta\lambda}{8\pi}\xi \int_r^{+\infty} dx \int_{-\infty}^x \exp\left(2 \left[V_\beta(x) - V_\beta(y)\right]\right) 
g_\beta(y)dy\,.
\end{align}
Noting that $g_\beta(r)\le 1$, we obtain the lower bound
 \begin{align}\label{first-ineq}
1- \xi \frac{\beta\lambda}{8\pi} t_\beta(r) \le   g_\beta(r) \,. 
\end{align}

 Using the asymptotic \eqref{equiv-t-beta}, we obtain a lower bound
 \begin{align*}
1- \xi \le   \liminf_{\beta\to 0} g_\beta(r_\beta)  \,. 
\end{align*}
Plugging \eqref{first-ineq} into \eqref{fixed-point-rescaled}, we get an upper bound 
\begin{align*}
g_\beta(r)  \le  1 - \xi \frac{\beta\lambda}{8\pi} t_\beta(r) +2  \xi^2 (\frac{\beta\lambda}{8\pi}) ^2 
\int_r^{+\infty} dx \int_{-\infty}^x \exp\left(2 \left[V_\beta(x) - V_\beta(y)\right]\right) t_\beta(y) dy\,. 
\end{align*}
After the derivation of the asymptotic quadruple integral  (similar to the one of \eqref{equiv-t-beta}), we 
obtain 
\begin{align*}
\limsup_{\beta\to 0} g_\beta(r_\beta) \le 1-\xi+ \xi^2\,. 
\end{align*}
Iterating this argument (the multiple integrals are always such that the integration ranges permit to 
catch the maximum value of $V_\beta(x)-V_\beta(y)$ inside the exponential as in \eqref{equiv-t-beta}), we get the result. 
\qed

\subsection{Proof of the estimates for a single diffusion}\label{proof-auxiliary-results-1}
Recall the definition of $R_{\lambda}:= \log(\tan(\alpha_{\lambda}/4))$ and its differential equation \eqref{eq-R}:
\begin{align*}
dR_\lambda= \frac{1}{2} \left(\lambda \frac{\beta}{4} e^{-\frac{\beta}{4}t} \cosh(R_\lambda) +  \tanh(R_\lambda) \right) dt + dB_t\,, 
\quad R_\lambda(0)=-\infty\,,
\end{align*}
and denote by $T_r(t)$ the first passage time to $r\in \R\cup \{+\infty\}$ after time $t$ i.e. \begin{equation*}T_r(t):=\inf\{s\geq t: R_\lambda(s) = r\}.\end{equation*} 
In this sub-section, we first study the diffusion $R_{\lambda}$ and then 
translate the estimates to the diffusion $\alpha_{\lambda}$.
We will denote by $\P_{r_0,t}$ the law of the diffusion $R_\lambda$ starting from position $r_0$ at time $t$. To simplify notations, we omit the subscript $t$ in $\P_{r_0,t}$ 
 if $t=0$ and write $T_r$ instead of $T_r(t)$ if it appears inside the probability $\P_{r_0,t}$.  
 We will also denote by $P_x$ the law of a Brownian motion starting from $x$.

Our first lemma \ref{fast-explosion} shows that when the diffusion $R_{\lambda}$ is outside of the well the probability that it explodes in a short time tends to $1$.
\begin{lemma}\label{fast-explosion}
Let $0 <\varepsilon<1$. 
Then,
 there exists a constant $c >0$ such that for all $\beta>0$ small enough,  
\begin{align*}
\P_{\varepsilon \log\frac{1}{\beta} } \left[ T_{+\infty} < 9    \log \frac{1}{\beta}  \right]  \geq 1- \beta^c  \,. 
\end{align*}
\end{lemma}
Note that it immediately gives the analogous result for $\alpha_{\lambda}$ of Lemma \ref{fast-explosion-diff}.

\noindent {\it Proof of Lemma \ref{fast-explosion}}.
We have 
\begin{align*}
\P_{\varepsilon \log\frac{1}{\beta} } \left[ T_{+\infty} < 9    \log \frac{1}{\beta}  \right] \ge 
\P_{\varepsilon \log\frac{1}{\beta} } \left[ T_{2 \log\frac{1}{\beta}} < 8    \log \frac{1}{\beta} \wedge T_{\frac{\varepsilon}{2}\log\frac{1}{\beta}}  \right] \P_{2 \log\frac{1}{\beta},8\log\frac{1}{\beta}} \left[T_{+\infty} < \log\frac{1}{\beta}  \right]\,. 
\end{align*}
If $\frac{\varepsilon}{2} \ln \frac{1}{\beta} < r < 2\log \frac{1}{\beta}$, the drift term $-V_\beta'(r) \ge \frac{1}{2} \frac{e^{r}-e^{-r}}{e^r + e^{-r}} \ge \frac{1}{4}$ for $\beta$ small enough. Thus, for $t< T_{\frac{\varepsilon}{2}\log\frac{1}{\beta}}$, 
we have $R_\lambda(t) \ge B_t + t/4$. 
Therefore, 
\begin{align*}
\P_{\varepsilon \log\frac{1}{\beta} }& \left[ T_{2 \log\frac{1}{\beta}} < 8    \log \frac{1}{\beta} \wedge T_{\frac{\varepsilon}{2}\log\frac{1}{\beta}}  \right] 
\ge P_{\varepsilon \log\frac{1}{\beta} } \left[ \inf_{0 \le t \le 8\log\frac{1}{\beta}} B_t > \frac{\varepsilon}{2}\log\frac{1}{\beta} \right]\,.    
\end{align*}
This latter probability is easily computed with the reflection principle for Brownian motion. 
\begin{align*}
P_{\varepsilon \log\frac{1}{\beta} } &\left[ \inf_{0 \le t \le 8\log\frac{1}{\beta}} B_t > \frac{\varepsilon}{2}\log\frac{1}{\beta} \right]
= P_0\left[\sup_{0 \le t \le 8\log\frac{1}{\beta}} B_t \le \frac{\varepsilon}{2}\log\frac{1}{\beta}\right] \\ 
&=P_0 \left[ |B(1)| \le \frac{\varepsilon}{4\sqrt{2}} \sqrt{\log\frac{1}{\beta}} \right]
 = 1- O(\exp(-c\log\frac{1}{\beta}))\,.    
\end{align*}

We now consider the probability 
\begin{align*}
\P_{2 \log\frac{1}{\beta},8\log\frac{1}{\beta}} \left[T_{+\infty} < \log\frac{1}{\beta} \right]
\end{align*}
We define $G_\lambda(t):=R_\lambda(t)-B_t$ such that $G_\lambda(0)= 2\log\frac{1}{\beta}$. The function $G_\lambda$ satisfies the following random 
{\it ordinary} differential equation 
\begin{align*}
G_\lambda'(t)= \frac{1}{2}\left(\frac{\lambda\beta}{8} e^{-\frac{\beta}{4} t }\frac{e^{G_\lambda(t)+B_t} + e^{-G_\lambda(t)-B_t} }{2} 
+ \frac{e^{G_\lambda(t)+B_t} - e^{-G_\lambda(t)-B_t}}{e^{G_\lambda(t)+B_t} + e^{-G_\lambda(t)-B_t}} \right)\,. 
\end{align*} 
On the event 
\begin{align*}
\mathcal{E}_\beta:=\left\{\sup_{0\le t\le \log\frac{1}{\beta}} | B_t |  \le \frac{1}{2} \log \frac{1}{\beta} \right\}\,,
\end{align*}
which occurs with probability bigger than $1-\exp(-c\log(\frac{1}{\beta}))$ where $c>0$ is a constant independent of $\beta$, 
we can easily check that, for $t\in [0, \log\frac{1}{\beta}]$ and $\beta>0$ small enough,  
\begin{align*}
G_\lambda'(t)\ge \frac{\lambda\beta^{3/2}}{64} e^{G_\lambda(t)} - \frac{1}{2}\,. 
\end{align*}
This leads us to study the Cauchy problem 
\begin{align*}
H'_\lambda(t) =  \frac{\lambda\beta^{3/2}}{64} e^{H_\lambda(t)} - \frac{1}{2}, \quad H_\lambda(0)= 2\log\frac{1}{\beta}\,. 
\end{align*}
This ordinary differential equation can be solved explicitly. We find 
\begin{align*}
e^{-H_\lambda(t)} e^{-\frac{t}{2}} = \beta^{2} - \frac{\lambda\beta^{3/2}}{32} (1- e^{-\frac{t}{2}} ) \,. 
\end{align*}
It is easy to see that the exploding time of $H_\lambda$ is of order $\sqrt{\beta}$ as $\beta\to 0$. 
On the event $\mathcal{E}_\beta$, we have $H_\lambda(t) \le G_\lambda(t)$ for all $t\ge 0$. 

We can conclude that 
\begin{align*}
\P_{2 \log\frac{1}{\beta},8\log\frac{1}{\beta}} \left[T_{+\infty} < \log\frac{1}{\beta}  \right] \ge 
P_0\left[\mathcal{E}_\beta\right] = 1-\exp(-c\log(\frac{1}{\beta}))\,. 
\end{align*}

\qed

\begin{lemma}\label{lower-bound-proba-fast-explosion}
For all $\beta>0$ small enough,
\begin{align}\label{lower-bound-proba-f-e}
\P_{{-\frac{1}{4} \log\frac{1}{\beta}} } \left[T_{+\infty} < 10 \log\frac{1}{\beta}\right] \geq \sqrt{\beta}\,. 
\end{align}
\end{lemma}

\noindent {\it Proof of Lemma \ref{lower-bound-proba-fast-explosion}}.
The probability under consideration can be bounded below by 
\begin{align*}
\P_{-\frac{1}{4} \log\frac{1}{\beta}} \left[ T_{\frac{1}{4} \log\frac{1}{\beta}} <   \log \frac{1}{\beta}  \right]
\P_{\frac{1}{4} \log\frac{1}{\beta} ,\log\frac{1}{\beta}} \left[ T_{+\infty} < 9    \log \frac{1}{\beta} \right]  \,.
\end{align*}
For the first probability, we note that $R_\lambda(t) \geq -\frac{t}{2}+ B_t$ and thus 
\begin{align*}
\P_{-\frac{1}{4} \log\frac{1}{\beta}} \left[ T_{\frac{1}{4} \log\frac{1}{\beta}} <   \log \frac{1}{\beta}  \right] 
&\geq P_{0}    \left[ B( \log\frac{1}{\beta}) - \frac{1}{2} \log \frac{1}{\beta} \ge \frac{1}{2} \log \frac{1}{\beta} \right]\\
&= P_{0} \left[B_1 \geq  \sqrt{ \log\frac{1}{\beta} }\right]\sim  \frac{\sqrt{\beta}}{\sqrt{2\pi\log\frac{1}{\beta}}} \,. 
\end{align*}

From Lemma \ref{fast-explosion}, the second probability is bigger than $1-\beta^c$ for some positive constant $c>0$. 
\qed

On the large scale-time of the order $1/\beta$, we prove that the time spent by $R_{\lambda}$ in the well of the potential $V$ (near $0$ modulo $2\pi$ for $\alpha_{\lambda}$) is negligible, i.e. small compared to the typical time $1/\beta$ between two jumps. 
This is Lemma \ref{time-spent-bad-region}:
\begin{lemma}\label{time-spent-bad-region}
Let $t>0$ and 
\begin{align*}
\Xi_\beta(t) :=\E_{-\infty}\left[ \int_{0}^{\frac{8\pi }{\beta}t} 1_{\{R_\lambda (u) \geq - \frac{1}{4} \log\frac{1}{\beta} \}} du\right]\,. 
\end{align*}
Then, there exists $C>0$ independent of $\beta$ such that 
\begin{align*}
 \Xi_\beta(t)  \leq \frac{C}{\sqrt{\beta}} \log \frac{1}{\beta} \,. 
\end{align*}
\end{lemma}
Lemma \ref{time-spent-bad-region} can be translated for $\alpha_\lambda$, it gives Lemma \ref{expected-time-bad-region-alpha}.

{\it Proof of Lemma \ref{time-spent-bad-region}.}
The key estimate is the lower-bound \eqref{lower-bound-proba-f-e} given by Lemma \ref{lower-bound-proba-fast-explosion} 
on the probability for the diffusion $R_\lambda$ starting from the position 
$ -\frac{1}{4} \log\frac{1}{\beta}$ to blow up in a short time. The idea is then to relate  the time spent by the diffusion 
$R_\lambda$ above the level  $ -\frac{1}{4} \log\frac{1}{\beta}$ with the number of explosions in the interval 
$[0,\frac{8\pi t}{\beta}]$ which is of order $O(1)$ (the typical time between two explosions is of order $1/\beta$ 
from Lemma \ref{conv-exp}). 
For any $u\in \R$, we have 
\begin{align*}
&\P_{-\infty}\left[R_\lambda (u) \geq - \frac{1}{4} \log\frac{1}{\beta}\right]\\ & \leq 
\P_{-\infty}\left[R_\lambda (u) \geq - \frac{1}{4} \log\frac{1}{\beta}, T_{+\infty} <  10 \log\frac{1}{\beta} \right] + 
\P_{-\infty}\left[R_\lambda (u) \geq - \frac{1}{4} \log\frac{1}{\beta}, T_{+\infty} \ge 10 \log\frac{1}{\beta} \right]
\\& \leq  \P_{-\infty}\left[R_\lambda (u) \geq - \frac{1}{4} \log\frac{1}{\beta}, T_{+\infty}  < 10 \log\frac{1}{\beta} \right] +
\left(1-\sqrt{\beta}\right)\,  \P_{-\infty}\left[R_\lambda (u) \geq - \frac{1}{4} \log\frac{1}{\beta}\right]\,,
\end{align*}
where we have used the simple Markov property, the inequality \eqref{lower-bound-proba-f-e} as well as the increasing property 
in the third line. 
We deduce that 
\begin{align}
&\P_{-\infty}\left[R_\lambda (u) \geq - \frac{1}{4} \log \frac{1}{\beta}\right] \leq \frac{1}{\sqrt{\beta}} \, 
 \P_{-\infty}\left[R_\lambda (u) \geq - \frac{1}{4} \log\frac{1}{\beta}, T_{+\infty}  < 10 \log\frac{1}{\beta} \right] \notag \\
  & \leq \frac{1}{\sqrt{\beta}}  \, \P_{-\infty}\left[\mbox{The interval $[u,u+ 10 \log\frac{1}{\beta} ]$ contains at least one explosion}\right] 
  \label{control-at-time-u} \,. 
\end{align}
Denoting by $k=\mu_\lambda^\beta[0,t]$ the total number of explosions of the diffusion $R_\lambda$ in the interval 
 $[0,\frac{8\pi t}{\beta}]$ and by $0 < \zeta_1 < \cdots < \zeta_k <\frac{8\pi t}{\beta}$ 
 the explosion times, we easily see that almost surely 
 \begin{align*}
 \int_{0}^{\frac{8\pi t}{\beta}}  1_{\{\exists i : \zeta_i \in [u,u+ 10 \log\frac{1}{\beta} ] \}} du \le 10 \log\frac{1}{\beta}\, \left(1+ \mu_\lambda^\beta[0,t]\right) \,. 
\end{align*} 
Using this inequality and integrating \eqref{control-at-time-u} with respect to $u$, we finally obtain    
\begin{align*}
\Xi_\beta(t)  \leq\left(1+ \E[\mu_\lambda^\beta[0,t]] \right) \, 
 \frac{10}{\sqrt{\beta}}  \log\frac{1}{\beta} \, \leq  \left(1+\frac{\lambda}{2\pi}\right) \frac{10}{\sqrt{\beta}}   \log\frac{1}{\beta} \,. 
\end{align*}
\qed

\subsection{Proof of Lemma \ref{def-poissons}}\label{appendix-def-poissons}


Denote $\mathcal{C}_0(\R_+)$ the space of continuous functions $f:\R_+\to \R_+$ such that $f(t)\to 0$ as $t\to+ \infty$. 
To prove that $\mathcal{P}_\lambda$ is indeed a $(\mathcal{F}_t)$-Poisson process with the correct intensity, 
it suffices to check that its Laplace functional 
satisfies 
\begin{align}\label{lp-poisson}
\E\left[\exp(- \int_s^t f(u) \mathcal{P}_\lambda(du)  ) \big | \, \mathcal{F}_s \right] =
\exp\left(- \int_{s}^t (1-e^{-f(u)}) \lambda e^{-2 \pi u} du \right)\,.
\end{align} 
We have to check that the filtration  $(\mathcal{F}_t)$ does not contain too much information compared to the natural filtration 
denoted $(\mathcal{F}_t^\lambda)$ associated 
to the process $\mathcal{P}_\lambda$ only, for \eqref{lp-poisson} to remain valid. 

We denote by $(\mathcal{G}_t)$ the filtration associated to the complex Brownian motion $(Z_t)$ which drives the processes 
$\alpha_\lambda, \alpha_{\lambda'}$  and $\alpha_{\lambda''}$ according to \eqref{eq-alpha}.  

The conclusion of Theorem \ref{stopping-times-pp} is equivalent to the convergence of the Laplace functional 
of the point process $\mu_\lambda^\beta$ to the Laplace functional of the Poisson process $\mathcal{P}_\lambda$ (see Proposition 11.1.VIII (ii) of \cite{daley-vere-jones}). Let us fix $\eps >0$. We can deduce that
for any $\varphi\in \mathcal{C}_0(\R_+)$ and any $0 <s<t$, almost-surely: 
\begin{align}\label{cond-conv}
\E\left[\exp(- \int_{s+\eps}^t \varphi(u) \mu_\lambda^\beta(du)  ) \big | \, \mathcal{G}_{\frac{8 \pi }{\beta}s} \right] \rightarrow_{\beta\to 0} 
\exp\left(- \int_{s+ \eps}^t \left(1-e^{-\varphi(u)}\right) \lambda e^{-2 \pi u} du \right)\,. 
\end{align} 
The main point is that the sigma field $\mathcal{G}_{\frac{8 \pi }{\beta}s}$ already contains the information 
on the two diffusions $(\alpha_\lambda(\frac{8 \pi }{\beta} u))_{0\le u \le s}, (\alpha_{\lambda'}(\frac{8 \pi }{\beta} u))_{0\le u \le s}$ and $(\alpha_{\lambda''}(\frac{8 \pi }{\beta} u))_{0\le u \le s}$ up to time $s$ as they are strong solutions to the stochastic differential equation system \eqref{eq-alpha}. Note also that we need to introduce a small gap between $t$ and $t+\eps$ so that \ref{cond-conv} is valid: it is indeed possible that the position of $\alpha_\lambda((\frac{8 \pi }{\beta}) s)$ induces a jump quickly after $s$. This issue can be circumvented by examining an independent random time $U$ distributed uniformly over the interval $[s,s+\eps]$. The position of $\alpha_\lambda(\frac{8 \pi }{\beta} U)$ belongs to an interval of the type $[2k\pi, 2k\pi + 4 \arctan(\beta^{1/4})]$ with probability going to $1$ when $\beta \to 0$ (this is Lemma \ref{expected-time-bad-region-alpha}). A straightforward adaptation of the proof of Theorem \ref{stopping-times-pp} then shows that the conditional law of the measure $\mu_\lambda^{\beta}$ over $[s+\eps,\infty)$ converges to the Poisson measure of intensity $\lambda\exp(- 2\pi t)dt$ and \eqref{cond-conv} holds.

From the joint convergence of the triplet \eqref{joint-conv} along the sequence $(\beta_k)$, we can deduce for $\varphi, f,g,h\in \mathcal{C}_0(\R_+)$, 
\begin{align}
\lim_{k\to \infty} \E\left[\exp\left(-\int_{s+\eps}^t \varphi(u) \mu_\lambda^{\beta_k}(du)\right) \exp\left(-\int_0^s ( f(u) \mu_\lambda^{\beta_k}(du) 
+g(u)\mu_{\lambda'}^{\beta_k}(du) + h(u) \mu_{\lambda''- \lambda'}^{\beta_k}(du)) \right) \right] \notag \\ =
\E\left[\exp\left(-\int_{s+\eps}^t \varphi(u) \mathcal{P}_\lambda(du)\right) \exp\left(-\int_0^s ( f(u) \mathcal{P}_\lambda(du) 
+ g(u)\mathcal{P}_{\lambda'}(du) + h(u)\mathcal{P}_{\lambda''-\lambda'}(du)  ) \right)\right] \label{with-joint-conv-mu}  \,. 
\end{align}
On the other hand, using in addition \eqref{cond-conv}, we can check that  
\begin{align}
&\lim_{k\to \infty} \E\left[\exp\left(-\int_{s+\eps}^t \varphi(u) \mu_\lambda^{\beta_k}(du)\right) \exp\left(-\int_0^s ( f(u) \mu_\lambda^{\beta_k}(du) 
+g(u)\mu_{\lambda'}^{\beta_k}(du) + h(u) \mu_{\lambda''- \lambda'}^{\beta_k}(du)) \right) \right] \notag \\ 
&=  \exp\left(- \int_{s+\eps}^t \left(1-e^{-\varphi(u)}\right) \lambda e^{-2 \pi u} du \right) 
\E\left[\exp\left(-\int_0^s ( f(u) \mathcal{P}_\lambda(du) 
+ g(u)\mathcal{P}_{\lambda'}(du) + h(u)\mathcal{P}_{\lambda''-\lambda'}(du)  ) \right)\right] \label{with-joint-conv-cond}   \,. 
\end{align}
Gathering \eqref{with-joint-conv-mu} and \eqref{with-joint-conv-cond} and taking the limit $\eps \to 0$ (using that a.s. $\mathcal{P}_\lambda$ does not jump on $s$), we obtain \eqref{lp-poisson}.

\qed

\end{document}